\documentclass[12pt, reqno]{amsart}
\usepackage[utf8]{inputenc}
\usepackage[T1]{fontenc}
\usepackage{lmodern}
\usepackage{amsmath, amssymb, amsthm}
\usepackage{geometry}
\usepackage{graphicx}
\usepackage{hyperref}
\usepackage{xcolor}
\hypersetup{
    colorlinks=true,
    linkcolor=blue,
    citecolor=blue,
    urlcolor=blue,
}

\newtheorem{theorem}{Theorem}[section]
\newtheorem{lemma}[theorem]{Lemma}

\theoremstyle{definition}
\newtheorem{definition}[theorem]{Definition}
\theoremstyle{remark}
\newtheorem{remark}[theorem]{Remark}

\begin{document}

\title[Connected Lie Groups with a Left-Invariant Intrinsic Chronometer]{Connected Lie Groups with a Left-Invariant Intrinsic Chronometer}

\author{Ramazan G. Beden}
\address{Department of General Control Problems, Faculty of Mechanics and Mathematics, Lomonosov Moscow State University, Moscow 119991, Russia}
\email{ramazan.gediev@math.msu.ru}

\begin{abstract}
In this paper, we investigate left-invariant intrinsic chronometers (time functions) on connected finite-dimensional Lie groups. We formulate and prove sub-Lorentzian analogs of Berestovskii's classical direct and inverse theorems. A one-to-one correspondence is established between chronometers satisfying a proposed system of axioms (which includes upper semicontinuity and the condition that the interiors of chronological diamonds form a base for the topology) and left-invariant sub-Lorentzian structures satisfying a chronological condition. It is shown that an axiomatic chronometer uniquely determines a closed pointed convex cone in the Lie algebra, whose linear span satisfies the Chow--Rashevskii condition, and an upper semicontinuous, concave, positively homogeneous anti-norm on it. Conversely, we prove that any such algebraic structure satisfying a local chronological condition induces an intrinsic chronometer with the required properties. The proofs rely on local estimates in exponential coordinates, the Hahn--Banach separation theorem, and properties of weak limits in optimal control problems.
\end{abstract}

\maketitle

\section{Introduction}

The connection between the axiomatic geometry of metric spaces and smooth differential geometry is one of the fundamental themes in modern mathematics. A classical breakthrough in this area was the resolution of Hilbert's fifth problem (Gleason, Montgomery, Zippin), which showed that every locally compact topological group that is a topological manifold admits a unique analytic Lie group structure. 
Building on this foundation, in 1988, V.\,N.~Berestovskii proved a remarkable inverse theorem in his seminal paper~\cite{Berestovskii1988}: every left-invariant intrinsic metric on a connected Lie group $G$ that induces the manifold topology is generated by some left-invariant sub-Riemannian structure (a Carnot--Carathéodory--Finsler structure). In other words, the metric completely determines the generating distribution (satisfying the Chow--Rashevskii condition) and the norm on it.

In recent decades, sub-Lorentzian geometry has undergone active development (see the works of M.~Grochowski~\cite{Grochowski2009}, E.~Minguzzi~\cite{Minguzzi2019}, V.\,N.~Berestovskii and I.\,A.~Zubareva~\cite{BerestovskiiZubareva2015, BerZub2025}, Yu.\,L.~Sachkov and E.\,F.~Sachkova~\cite{Sachkov2023}, and L.\,V.~Lokutsievskiy and A.\,V.~Podobryaev~\cite{Lokutsievskiy2024}). Unlike the sub-Riemannian case, instead of a distribution of norms, one is given a field of convex cones $C^+$ and concave anti-norms $\nu$. In parallel, in 2018, M.~Kunzinger and C.~Sämann~\cite{KunzingerSaemann2018} laid the foundations of synthetic Lorentzian geometry (Lorentzian length spaces).

Despite the active study of the direct problem (investigating reachability properties, geodesics, and global hyperbolicity for given sub-Lorentzian structures), the inverse problem for sub-Lorentzian Lie groups has remained open. The difficulties in solving it stem from the non-compactness of light cones, the lack of an analog to the Hopf--Rinow theorem for bounded sets, and the fact that the chronometer $\tau$ (the Lorentzian time function) is upper semicontinuous rather than continuous. In a topological context, the connection between left-invariant orders and semigroups was studied by V.\,N.~Berestovskii and V.\,M.~Gichev~\cite{BerestovskiiGichev2000}; however, the transition to a left-invariant sub-Lorentzian structure along with a proof of the Chow--Rashevskii condition was not achieved.

The purpose of the present paper is to close this gap and establish a strict equivalence between axiomatized left-invariant Lorentzian spaces and left-invariant sub-Lorentzian structures on connected Lie groups. 
In this paper, the author formulates a system of axioms for an intrinsic left-invariant chronometer $\tau$. 
In the main (inverse) Theorem~\ref{thm:main}, it is proved that such a chronometer uniquely determines a nontrivial, closed, pointed convex cone $C_0^+ \subset L = T_eG$ and an anti-norm $\nu_0$ on it, with the linear span of the cone $C_0^+$ generating the entire Lie algebra $L$. 
The direct Theorem~\ref{thm:direct} establishes the converse: any left-invariant sub-Lorentzian structure $(G, C_0^+, \nu_0)$ satisfying the proposed chronological condition generates an intrinsic chronometer that satisfies the stated system of axioms.

\vspace{1ex}
\noindent\textbf{Geometric idea of the proof.} 
Locally, a chronological diamond behaves like a convex cone. Instead of directly analyzing a sequence of shrinking diamonds, whose topological boundary may have a complicated structure, our approach is based on intersecting the closed chronological future with a sequence of shrinking Euclidean balls and subsequently dilating them in the tangent space (the Lie algebra). In the limit, this procedure defines a tangent cone. The pointedness of this cone prevents admissible curves from ``winding'' infinitely, which would accumulate infinite sub-Riemannian length without a strictly positive time increment. Geometrically, this property guarantees the fulfillment of the Chow--Rashevskii condition. Once the anti-norm is defined on the obtained cone, it remains to prove that the chronometer $\tau$ is exactly recovered by integrating this anti-norm along admissible curves.

\section{Local Coordinates and Estimates in a Lie Group (Lemmas 1--6 of V.\,N.~Berestovskii)}

In this section, we present local estimates in exponential coordinates. They are based on the classical lemmas of V.\,N.~Berestovskii \cite{Berestovskii1988}; however, for the reader's convenience and the needs of subsequent exposition, we provide them here with complete proofs.

Let $G$ be a connected Lie group with identity $e$ and Lie algebra $L = T_e G$. 
Consider a chart $(U, \varphi)$ of the Lie group $G$, where $U$ is an open neighborhood of the identity $e$ in $G$, and $\varphi = (\exp|_V)^{-1} : U \to V \subset L$ is a diffeomorphism onto an open ball $V = U(0, \lambda) = \{x \in L : |x| < \lambda\}$ in $L$ with respect to some Euclidean norm $|\cdot|$ on $L$.

For elements $x, y \in V$ such that $\exp x \cdot \exp y \in U$, we set
\begin{equation}
x * y := \varphi(\exp x \cdot \exp y) = x + y + z(x, y), \label{eq:ber1}
\end{equation}
where $z(x, y)$ is an analytic vector function on its domain of definition in $V \times V$, such that $z(x, 0) = z(0, y) = 0$.

\begin{lemma}\label{lember:1}
For a suitable choice of the Euclidean norm $|\cdot|$ on the Lie algebra $L$ and a number $\lambda \in (0, 1)$, for all $x, y \in V = U(0, \lambda)$ such that $\exp x \cdot \exp y \in U$, the following estimate holds:
\begin{equation}
|z(x, y)| \le |x| \cdot |y| \qquad (\text{where } |z(x, y)| < |x| \cdot |y| \text{ for } x \neq 0,\, y \neq 0). \label{eq:ber2}
\end{equation}
\end{lemma}

\begin{proof}
Since $x * 0 = x$ and $0 * y = y$, the mapping $z(x, y)$ vanishes when $x = 0$ and when $y = 0$. Therefore, in the Taylor expansion of the coordinate functions of $z(x, y)$ in a neighborhood of $(0, 0)$, there are no terms depending solely on $x$ or solely on $y$. By Taylor's formula, in some neighborhood of zero $U(0, \lambda_0) \times U(0, \lambda_0)$, the estimate $|z(x, y)| \le C |x| \cdot |y|$ holds for some constant $C > 0$. 
Passing to a new Euclidean norm $|x|_{\text{new}} := c |x|$ with $c > C$, we obtain
\[
|z(x, y)|_{\text{new}} = c |z(x, y)| \le c C |x| \cdot |y| = \frac{C}{c} |x|_{\text{new}} \cdot |y|_{\text{new}} \le |x|_{\text{new}} \cdot |y|_{\text{new}},
\]
where the inequality is strict for $x \neq 0$ and $y \neq 0$. Decreasing $\lambda \in (0, 1)$ if necessary, we obtain the required assertion.
\end{proof}

\begin{lemma}\label{lember:2}
Under the conditions of Lemma \ref{lember:1}, for any $x, y \in V$ such that $\exp x \cdot \exp y \in U$, the following inequality holds:
\[
|x * y| \le |x| + |y| + |x| \cdot |y| \qquad (< |x| + |y| + |x| \cdot |y| \text{ for } x, y \neq 0).
\]
\end{lemma}

\begin{proof}
The statement follows directly from the representation \eqref{eq:ber1}, estimate \eqref{eq:ber2}, and the triangle inequality:
\[
|x * y| = |x + y + z(x, y)| \le |x| + |y| + |z(x, y)| \le |x| + |y| + |x| \cdot |y|. \qedhere
\]
\end{proof}

\begin{lemma}\label{lember:3}
If $x, y, x + y \in V$, $(-x) * (x + y) \in V$, and $|x| < 1/2$, then
\[
x + y = x * y',
\]
where the vector $y' \in V$ satisfies the estimate
\[
|y - y'| \le 2 |x| \cdot |y| \qquad (< 2 |x| \cdot |y| \text{ for } x, y \neq 0).
\]
\end{lemma}

\begin{proof}
Let $y' := (-x) * (x + y)$. Then $x * y' = x * ((-x) * (x + y)) = x + y$. By formula \eqref{eq:ber1}, we have
\[
x + y = x * y' = x + y' + z(x, y'),
\]
which yields $y - y' = z(x, y')$. Applying Lemma \ref{lember:1} and the triangle inequality, we find
\[
|y - y'| = |z(x, y')| \le |x| \cdot |y'| \le |x| \big(|y| + |y - y'|\big).
\]
Moving the term $|x| \cdot |y - y'|$ to the left-hand side and taking into account the condition $|x| < 1/2$ (that is, $1 - |x| > 1/2$), we obtain
\[
|y - y'| \le \frac{|x| \cdot |y|}{1 - |x|} \le 2 |x| \cdot |y|,
\]
where the inequality is strict for $x, y \neq 0$.
\end{proof}

\begin{lemma}\label{lember:4}
If the elements $x, y, y * x, (-x) * y * x$ belong to $V$ and $|x| < 1/2$, then
\[
y * x = x * y',
\]
where the vector $y' \in V$ satisfies the estimate
\[
|y - y'| \le 4 |x| \cdot |y| \qquad (< 4 |x| \cdot |y| \text{ for } x, y \neq 0).
\]
\end{lemma}

\begin{proof}
Let $y' := (-x) * y * x$, so that $y * x = x * y'$. Expanding both sides of the equality $y * x = x * y'$ using formula \eqref{eq:ber1}, we have
\[
y + x + z(y, x) = x + y' + z(x, y'),
\]
which yields
\[
y - y' = z(x, y') - z(y, x).
\]
Estimating the right-hand side using Lemma \ref{lember:1} and the triangle inequality, we obtain
\[
|y - y'| \le |z(x, y')| + |z(y, x)| \le |x| \cdot |y'| + |y| \cdot |x| \le |x| \big(2|y| + |y - y'|\big).
\]
Consequently, $(1 - |x|)|y - y'| \le 2 |x| \cdot |y|$. Since $|x| < 1/2$, we have $1 - |x| > 1/2$, which implies
\[
|y - y'| \le \frac{2 |x| \cdot |y|}{1 - |x|} \le 4 |x| \cdot |y|,
\]
where the inequality is strict for $x, y \neq 0$.
\end{proof}

\begin{lemma}\label{lember:5}
Let $x_1, \dots, x_j \in V$, $|x_i| \le \eta / j$ ($i = 1, \dots, j$) for some $\eta > 0$, and let $(1 + \eta / j)^j - 1 < \lambda$. Then the product $x_1 * \dots * x_j$ is defined in $V$ and satisfies the estimate
\[
|x_1 * \dots * x_j| \le \left(1 + \frac{\eta}{j}\right)^j - 1 < e^\eta - 1.
\]
\end{lemma}

\begin{proof}
We prove the inequality 
\[
|x_1 * \dots * x_m| \le \left(1 + \frac{\eta}{j}\right)^m - 1
\]
by induction on $m = 1, \dots, j$. For $m = 1$, we have $|x_1| \le \eta / j = (1 + \eta / j)^1 - 1$. Suppose the estimate holds for $m - 1 < j$. By Lemma \ref{lember:2}, for the product $(x_1 * \dots * x_{m-1}) * x_m$, we obtain:
\begin{align*}
|x_1 * \dots * x_m| &\le |x_1 * \dots * x_{m-1}| + |x_m| + |x_1 * \dots * x_{m-1}| \cdot |x_m| \\
&\le \left(\left(1 + \frac{\eta}{j}\right)^{m-1} - 1\right) + \frac{\eta}{j} + \left(\left(1 + \frac{\eta}{j}\right)^{m-1} - 1\right)\frac{\eta}{j} \\
&= \left(1 + \frac{\eta}{j}\right)^{m-1}\left(1 + \frac{\eta}{j}\right) - 1 = \left(1 + \frac{\eta}{j}\right)^m - 1.
\end{align*}
For $m = j$, we obtain $|x_1 * \dots * x_j| \le (1 + \eta / j)^j - 1 < e^\eta - 1$.
\end{proof}

\begin{lemma}\label{lember:6}
Let $x_i, y_i \in V$, $|x_i| \le \eta / j$, $|y_i| \le \varepsilon \eta / j$ ($i = 1, \dots, j$) for some $\eta > 0$, where $0 < \varepsilon \le \min\left(1, \frac{\lambda}{14\eta e^\eta}\right)$, and assume the condition $e^{(1 + \varepsilon)\eta} - 1 < \lambda < 1/2$ is satisfied. Then
\[
\prod_{i=1}^j (x_i + y_i) = \left(\prod_{i=1}^j x_i\right) * y,
\]
where the product is understood with respect to the operation $*$, and the remainder vector $y \in V$ satisfies the estimate
\[
|y| < 14 \varepsilon \eta e^\eta.
\]
\end{lemma}

\begin{proof}
1. First, note that for each $i = 1, \dots, j$, the vector $(-x_i)*(x_i + y_i)$ is well-defined and lies in $V$. Indeed, by Lemma \ref{lember:2}, we have:
\[
|(-x_i)*(x_i + y_i)| \le |-x_i| + |x_i + y_i| + |-x_i||x_i + y_i| \le \frac{\eta}{j} + \frac{(1+\varepsilon)\eta}{j} + \frac{(1+\varepsilon)\eta^2}{j^2} \le \left(1 + \frac{(1+\varepsilon)\eta}{j}\right)^2 - 1.
\]
For $j \ge 2$, the right-hand side does not exceed $\big(1 + (1+\varepsilon)\eta/j\big)^j - 1 < e^{(1+\varepsilon)\eta} - 1 < \lambda$ (for $j = 1$, the vector $y_1'$ in the representation $x_1 + y_1 = x_1 * y_1'$ is also well-defined and lies in $V$, since $|x_1| + |x_1 + y_1| + |x_1||x_1 + y_1| \le e^\eta e^{(1+\varepsilon)\eta} - 1$).
Since $|x_i| \le \eta/j < e^\eta - 1 < \lambda < 1/2$, by Lemma \ref{lember:3}, we represent the sum $x_i + y_i$ in the form
\[
x_i + y_i = x_i * y_i', \qquad \text{where } |y_i - y_i'| \le 2 |x_i| \cdot |y_i|.
\]
Hence, taking into account that $|x_i| < 1/2$, we obtain the norm estimate for the vector $y_i'$:
\[
|y_i'| \le |y_i| + |y_i - y_i'| \le (1 + 2|x_i|)|y_i| \le 2|y_i| \le \frac{2\varepsilon\eta}{j}.
\]

2. Using the associativity of the operation $*$, we move all factors $y_i'$ to the right side of the product $\prod_{i=1}^j (x_i * y_i')$. To do this, let $y_j'' := y_j'$, and for each $i = 1, \dots, j - 1$, define the vector $y_i''$ from the relation
\[
y_i' * X_i = X_i * y_i'', \qquad \text{where } X_i := x_{i+1} * \dots * x_j.
\]
Due to the conditions of the lemma, all intermediate products are guaranteed to lie in the ball $V$; therefore, the application of Lemma \ref{lember:4} is valid. Sequential rearrangement yields
\[
\prod_{i=1}^j (x_i + y_i) = (x_1 * y_1') * (x_2 * y_2') * \dots * (x_j * y_j') = \left(\prod_{i=1}^j x_i\right) * (y_1'' * y_2'' * \dots * y_j'').
\]

3. Let us estimate the norms of the vectors $y_i''$. By Lemma \ref{lember:5}, we have
\[
|X_i| = |x_{i+1} * \dots * x_j| \le \left(1 + \frac{\eta}{j}\right)^{j-i} - 1 \le \left(1 + \frac{\eta}{j}\right)^j - 1 < e^\eta - 1,
\]
where, by assumption, $|X_i| < \lambda < 1/2$. Applying Lemma \ref{lember:4} to the equality $y_i' * X_i = X_i * y_i''$, we find
\[
|y_i' - y_i''| \le 4 |X_i| \cdot |y_i'|,
\]
which yields
\[
|y_i''| \le (1 + 4|X_i|)|y_i'| \le \big(1 + 4(e^\eta - 1)\big)|y_i'| = (4e^\eta - 3)|y_i'| \le 4e^\eta |y_i'| \le \frac{8 \varepsilon \eta e^\eta}{j}.
\]

4. Denote $y := y_1'' * y_2'' * \dots * y_j''$ and let $\tilde{\eta} := 8 \varepsilon \eta e^\eta > 0$. By the condition of the lemma, $\varepsilon \le \frac{\lambda}{14\eta e^\eta}$, so $\tilde{\eta} \le \frac{8}{14}\lambda < 1$ and $(1 + \tilde{\eta}/j)^j - 1 < e^{\tilde{\eta}} - 1 \le (e - 1)\tilde{\eta} < 14\varepsilon \eta e^\eta \le \lambda$. Since $|y_i''| \le \tilde{\eta} / j$ for all $i = 1, \dots, j$, applying Lemma \ref{lember:5} to the product $y$ gives $y \in V$ and the estimate
\[
|y| \le \left(1 + \frac{\tilde{\eta}}{j}\right)^j - 1 < e^{\tilde{\eta}} - 1 \le (e - 1)\tilde{\eta} = 8(e - 1)\varepsilon \eta e^\eta < 14 \varepsilon \eta e^\eta. \qedhere
\]
\end{proof}

\section{Basic Definitions}

Let $(G, \tau)$ be a connected finite-dimensional Lie group with identity $e$ and a chronometer $\tau : G \times G \to [0, +\infty] \cup \{-\infty\}$ satisfying the following properties:

\begin{enumerate}
\item \textbf{Left-invariance:} for any $a, g, h \in G$, the following equality holds:
\[
\tau(ag, ah) = \tau(g, h).
\]

\item \textbf{Reflexivity and antisymmetry:} for any $g \in G$, $\tau(g, g) = 0$; if $g \neq h$ and $\tau(g, h) \ge 0$, then $\tau(h, g) = -\infty$.

\item \textbf{Reverse triangle inequality:} for any $g_1, g_2, g_3 \in G$ such that $\tau(g_1, g_2) \ge 0$ and $\tau(g_2, g_3) \ge 0$, the following inequality holds:
\[
\tau(g_1, g_3) \ge \tau(g_1, g_2) + \tau(g_2, g_3).
\]

\item \textbf{Local upper semicontinuity:} there exists an open neighborhood $U_0$ of the identity $e \in G$ such that for any $c \ge 0$, the set
\[
\{(g, h) \in U_0 \times U_0 : \tau(g, h) \ge c\}
\]
is closed in $U_0 \times U_0$.

\item \textbf{Consistency with the manifold topology:} for any $p, q \in G$, let
\[
I^+(p) := \{x \in G : \tau(p, x) > 0\}, \qquad I^-(q) := \{x \in G : \tau(x, q) > 0\},
\]
\[
I(p, q) := I^+(p) \cap I^-(q) = \{x \in G : \tau(p, x) > 0,\ \tau(x, q) > 0\}.
\]
The family of interiors of chronological diamonds $\{\operatorname{int} I(p, q)\}_{p, q \in G}$ forms a base for the topology of the group $G$, and for any point $x \in G$ and any neighborhood $U \subset G$ of $x$, there exist points $p, q \in U$ such that
\[
x \in \operatorname{int} I(p, q) \subset I(p, q) \subset U.
\]

\item \textbf{Intrinsic nature of the chronometer:} a continuous curve $\gamma : [a, b] \to G$ is called $\tau$-\textit{admissible} if for any $s_1, s_2 \in [a, b]$, $s_1 < s_2$ implies $\tau(\gamma(s_1), \gamma(s_2)) \ge 0$. The $\tau$-\textit{length} of such a curve is defined as the quantity
\[
L_\tau(\gamma) := \inf \sum_{i=1}^m \tau\big(\gamma(s_{i-1}), \gamma(s_i)\big),
\]
where the infimum is taken over all finite partitions $a = s_0 < s_1 < \dots < s_m = b$ of the interval $[a, b]$. For any distinct points $g, h \in G$ satisfying $\tau(g, h) \ge 0$, the set $\Omega_\tau(g, h)$ of all continuous $\tau$-admissible curves connecting $g$ to $h$ is nonempty, and the following equality holds:
\[
\tau(g, h) = \sup_{\gamma \in \Omega_\tau(g, h)} L_\tau(\gamma).
\]
\end{enumerate}

\begin{remark}
The sets
\[
J^+(p) := \{x \in G : \tau(p, x) \ge 0\}, \qquad J^-(q) := \{x \in G : \tau(x, q) \ge 0\}
\]
are called the closed future of the point $p$ and the closed past of the point $q$, respectively. The set
\[
J(p, q) := J^+(p) \cap J^-(q) = \{x \in G : \tau(p, x) \ge 0,\ \tau(x, q) \ge 0\}
\]
is called the closed chronological diamond of the points $p$ and $q$.
\end{remark}

\begin{remark}
Note that the reverse triangle inequality for the chronometer and its antisymmetry imply the absence of nontrivial closed $\tau$-admissible curves.
\end{remark}

\begin{definition}
A left-invariant sub-Lorentzian structure on a Lie group $G$ with Lie algebra $L = T_e G$ is a triple $(G, C_0^+, \nu_0)$, where $C_0^+ \subset L$ is a convex, nontrivial ($C_0^+ \neq \{0\}$), closed, pointed ($C_0^+ \cap (-C_0^+) = \{0\}$) cone whose linear span $L_0 := C_0^+ - C_0^+$ generates the entire Lie algebra $L$, and $\nu_0 : L \to [0, +\infty) \cup \{-\infty\}$ is an upper semicontinuous function (anti-norm) satisfying the following properties:
\begin{enumerate}
\item \textbf{Domain of nonnegativity and nontriviality:} $\nu_0(u) \ge 0$ for all $u \in C_0^+$, and $\nu_0(u) = -\infty$ for all $u \in L \setminus C_0^+$; moreover, $\nu_0(u) > 0$ on the relative interior $\operatorname{int}_{L_0}(C_0^+)$ (equivalently, $\nu_0|_{C_0^+} \not\equiv 0$);
\item \textbf{Positive homogeneity:} for any $u \in L$ and $\lambda > 0$, the following equality holds:
\[
\nu_0(\lambda u) = \lambda \nu_0(u);
\]
\item \textbf{Concavity:} for any $u, v \in L$, the following inequality holds:
\[
\nu_0(u + v) \ge \nu_0(u) + \nu_0(v).
\]
\end{enumerate}
\end{definition}

Let $l_g : G \to G$, $l_g(h) := gh$, denote the left translation operator by an element $g \in G$, and let $(dl_g)_e : L \to T_g G$ be its differential at the identity $e$. Left translations of the cone $C_0^+$ and the anti-norm $\nu_0$ define a left-invariant cone field $C_g^+ := (dl_g)_e(C_0^+) \subset T_g G$ and a left-invariant anti-norm field $\nu_g(\xi) := \nu_0\big((dl_{g^{-1}})_g(\xi)\big)$, $\xi \in T_g G$, on the Lie group $G$.

An absolutely continuous curve $\gamma : [0, 1] \to G$ is said to be $C^+$-\textit{admissible} (or $C_0^+$-\textit{admissible}) if for almost all $s \in [0, 1]$ its tangent vector $\dot{\gamma}(s)$ belongs to the cone $C_{\gamma(s)}^+$, that is, its left-translated velocity vector to the identity
\[
\theta\big(\dot{\gamma}(s)\big) := \big(dl_{\gamma(s)^{-1}}\big)_{\gamma(s)}\big(\dot{\gamma}(s)\big) \in L
\]
lies in the cone $C_0^+$ almost everywhere on $[0, 1]$.

The \textit{sub-Lorentzian length} of a $C^+$-admissible curve $\gamma : [0, 1] \to G$ is defined as the quantity
\[
\ell(\gamma) := \int_0^1 \nu_{\gamma(s)}\big(\dot{\gamma}(s)\big)\,ds = \int_0^1 \nu_0\big(\theta(\dot{\gamma}(s))\big)\,ds.
\]
By the positive homogeneity of the anti-norm $\nu_0$, the sub-Lorentzian length $\ell(\gamma)$ is independent of the choice of a monotonically nondecreasing absolutely continuous parametrization of the curve $\gamma$.

For any points $g, h \in G$, we denote by $\Omega_{C^+}(g, h)$ the set of all $C^+$-admissible curves $\gamma : [0, 1] \to G$ connecting $g$ to $h$ ($\gamma(0) = g$, $\gamma(1) = h$). The \textit{sub-Lorentzian chronometer} generated by the structure $(G, C_0^+, \nu_0)$ is the function $\tau_c : G \times G \to [0, +\infty] \cup \{-\infty\}$ defined by the equality
\[
\tau_c(g, h) := \begin{cases}
\sup\limits_{\gamma \in \Omega_{C^+}(g, h)} \ell(\gamma), & \text{if } \Omega_{C^+}(g, h) \neq \varnothing, \\[2ex]
-\infty, & \text{if } \Omega_{C^+}(g, h) = \varnothing.
\end{cases}
\]

\section{Main Theorem}

We state the main result of this paper.

\begin{theorem}\label{thm:main}
Let $(G, \tau)$ be a connected finite-dimensional Lie group with a chronometer $\tau$ satisfying properties 1--6. Then there exist a unique closed pointed convex cone $C_0^+$ in the Lie algebra $L$ and a unique anti-norm $\nu_0$ on the cone $C_0^+$ such that the triple $(G, C_0^+, \nu_0)$ forms a left-invariant sub-Lorentzian structure on the Lie group $G$ generating the chronometer $\tau$:
\[
\tau(g, h) = \tau_c(g, h) \qquad \text{for all } g, h \in G.
\]
Moreover, the cone $C_0^+$ and the anti-norm $\nu_0$ on $C_0^+$ are uniquely determined by the relations
\[
C_0^+ = \left\{u \in L : \tau\big(e, \exp(tu)\big) \ge 0 \text{ for all sufficiently small } t > 0\right\},
\]
\[
\nu_0(u) = \lim_{t \to 0^+} \frac{\tau\big(e, \exp(tu)\big)}{t} = \inf_{t > 0} \frac{\tau\big(e, \exp(tu)\big)}{t} \qquad (u \in C_0^+).
\]
\end{theorem}

\begin{proof}
Fix a Euclidean norm $|\cdot|$ on $L$ and a number $\lambda \in (0, 1)$ such that the exponential map is a diffeomorphism from the open ball $V = U(0, \lambda) = \{x \in L : |x| < \lambda\}$ onto an open neighborhood $U \subset U_0$ of the identity in $G$, the chronometer $\tau$ is upper semicontinuous and finite on $U \times U$ ($\tau(g, h) < 1$ for all $g, h \in U$, which can be achieved by decreasing $\lambda$ due to $\tau(e, e) = 0$ and the closedness of the set $\{\tau \ge 1\}$ in $U_0 \times U_0$), and the group operation $x * y = \exp^{-1}(\exp(x)\exp(y)) = x + y + z(x, y)$ satisfies the estimate of Berestovskii's Lemma~\ref{lember:1}:
\[
|z(x, y)| \le |x| \cdot |y|
\]
for all $x, y \in V$ such that $x * y \in V$. For each $\rho \in (0, \lambda)$, denote by $B_\rho(0) = \{x \in L : |x| \le \rho\}$ the closed ball in $L$, and by $B^G_\rho(e) = \exp(B_\rho(0))$ its image under the exponential map.

We prove the following lemma.

\begin{lemma}\label{lem:1}
There exist points $p_0, q_0 \in U$ and a compact neighborhood $K \subset U$ of the identity $e$ such that
\[
e \in \operatorname{int} I(p_0, q_0) \subset I(p_0, q_0) \subset K;
\]
moreover, for any points $p, q \in I(p_0, q_0)$, the set $J(p, q)$ is a compact subset of $I(p_0, q_0) \subset K$. More generally, for any nonempty compact sets $K_1, K_2 \subset I(p_0, q_0)$, the set
\[
J(K_1, K_2) := \bigcup_{p \in K_1,\, q \in K_2} J(p, q)
\]
is compact and is contained in $I(p_0, q_0) \subset K$.
\end{lemma}

\begin{proof}[Proof of Lemma \ref{lem:1}]
Choose an arbitrary compact neighborhood $K \subset U$ of the identity $e$ (for example, $K = B_{\rho}^G(e)$ for some $\rho \in (0, \lambda)$). By property 5, for the point $e$ and its open neighborhood $\operatorname{int} K$, there exist points $p_0, q_0 \in \operatorname{int} K \subset U$ such that
\[
e \in \operatorname{int} I(p_0, q_0) \subset I(p_0, q_0) \subset \operatorname{int} K \subset K \subset U.
\]
Let $p, q \in I(p_0, q_0)$; then $\tau(p_0, p) > 0$ and $\tau(q, q_0) > 0$. For any point $x \in J(p, q)$, we have $\tau(p, x) \ge 0$ and $\tau(x, q) \ge 0$, which implies, by the reverse triangle inequality (property 3):
\[
\tau(p_0, x) \ge \tau(p_0, p) + \tau(p, x) \ge \tau(p_0, p) > 0, \qquad
\tau(x, q_0) \ge \tau(x, q) + \tau(q, q_0) \ge \tau(q, q_0) > 0.
\]
Therefore, $x \in I(p_0, q_0)$, that is, $J(p, q) \subset I(p_0, q_0) \subset K \subset U$.

Now let $K_1, K_2 \subset I(p_0, q_0)$ be nonempty compact sets. As proved above, $J(K_1, K_2) \subset I(p_0, q_0) \subset K$. We show that the set $J(K_1, K_2)$ is closed in the compact set $K$. Let a sequence $\{x_n\} \subset J(K_1, K_2)$ converge to a point $x \in K$. For each $n$, there exist points $p_n \in K_1$ and $q_n \in K_2$ such that $\tau(p_n, x_n) \ge 0$ and $\tau(x_n, q_n) \ge 0$. By the compactness of $K_1$ and $K_2$, passing to a subsequence, we may assume that $p_n \to p \in K_1$ and $q_n \to q \in K_2$ as $n \to \infty$. Since $K_1, K_2, K \subset U \subset U_0$, and by property 4 the set $\{(g, h) \in U \times U : \tau(g, h) \ge 0\}$ is closed in $U \times U$, taking the limit as $n \to \infty$ yields $\tau(p, x) \ge 0$ and $\tau(x, q) \ge 0$. Thus, $x \in J(p, q) \subset J(K_1, K_2)$, so $J(K_1, K_2)$ (and, in particular, $J(p, q)$ for $K_1 = \{p\}$ and $K_2 = \{q\}$) is compact.
\end{proof}

We prove the additivity property of the $\tau$-length functional $L_\tau$ as a function of the interval.

\begin{lemma}\label{lem:2}
For any continuous $\tau$-admissible curve $\gamma : [a, b] \to G$ and any point $c \in (a, b)$, the following equality holds:
\[
L_\tau\big(\gamma|_{[a, b]}\big) = L_\tau\big(\gamma|_{[a, c]}\big) + L_\tau\big(\gamma|_{[c, b]}\big).
\]
\end{lemma}

\begin{proof}[Proof of Lemma \ref{lem:2}]
The union of any partition of the interval $[a, c]$ and any partition of the interval $[c, b]$ forms a partition of the interval $[a, b]$ containing the point $c$. Therefore, by the definition of the infimum, we have $L_\tau(\gamma|_{[a, b]}) \le L_\tau(\gamma|_{[a, c]}) + L_\tau(\gamma|_{[c, b]})$. Conversely, if $P = \{a = s_0 < s_1 < \dots < s_m = b\}$ is an arbitrary partition of the interval $[a, b]$, then adding the point $c \in [s_{k-1}, s_k]$ to it yields, by the reverse triangle inequality (property 3), $\tau(\gamma(s_{k-1}), \gamma(s_k)) \ge \tau(\gamma(s_{k-1}), \gamma(c)) + \tau(\gamma(c), \gamma(s_k))$. Consequently, the sum over the partition $P$ is not less than the sum over the refined partition $P \cup \{c\}$, which in turn is not less than $L_\tau(\gamma|_{[a, c]}) + L_\tau(\gamma|_{[c, b]})$. Taking the infimum over all partitions $P$, we obtain the reverse inequality $L_\tau(\gamma|_{[a, b]}) \ge L_\tau(\gamma|_{[a, c]}) + L_\tau(\gamma|_{[c, b]})$.
\end{proof}

By property 5, for the point $e$ and its open neighborhood $\operatorname{int} I(p_0, q_0)$, we choose points $p, q \in I(p_0, q_0)$ such that
\[
e \in \operatorname{int} I(p, q) \subset I(p, q) \subset I(p_0, q_0).
\]
Then $\tau(e, q) > 0$, which implies $q \neq e$ by property 2, and by Lemma~\ref{lem:1}, the set $J(e, q)$ is compact. Since $\operatorname{int} I(p, q)$ is an open neighborhood of the identity $e$ and $q \neq e$, we choose $r > 0$ small enough so that the following conditions are satisfied:
\[
e^{2r} - 1 < \min\left(\lambda, \frac{1}{2}\right), \qquad B^G_r(e) = \exp(B_r(0)) \subset \operatorname{int} I(p, q) \subset I(p, q), \qquad q \notin B^G_r(e).
\]
In particular, for all $g \in B^G_r(e)$, the inequality $\tau(g, q) > 0$ holds.

Let us proceed to construct the cone. Consider the sets
\[
M_{\frac{r}{2^k}} := J(e, q) \cap B^G_{\frac{r}{2^k}}(e) = \left\{ g \in B^G_{\frac{r}{2^k}}(e) : \tau(e, g) \ge 0 \right\}, \qquad k \ge 1.
\]
Let $A_{\frac{r}{2^k}}$ denote $\exp^{-1}\big(M_{\frac{r}{2^k}}\big)$ and set
\[
C_k := 2^k A_{\frac{r}{2^k}} = \left\{ a \in B_r(0) : \tau\left(e, \exp\left( \frac{a}{2^k} \right)\right) \ge 0 \right\}.
\]
Since $J(e, q)$ and $B^G_{\frac{r}{2^k}}(e)$ are compact, and $\exp|_{B_r(0)}$ is a homeomorphism onto $B^G_r(e)$, the sets $M_{\frac{r}{2^k}}$, $A_{\frac{r}{2^k}}$, and $C_k$ are nonempty compact sets (since $\tau(e, e) = 0$, the vector $0$ belongs to all $C_k$).

We show that $C_{k+1} \subset C_k$ for all $k \ge 1$. Let $x \in C_{k+1}$; then $|x| \le r$ and $\tau\left(e, \exp\left(\frac{x}{2^{k+1}}\right)\right) \ge 0$. By left-invariance (property 1) and the reverse triangle inequality (property 3), we obtain:
\[
\tau\left(e, \exp\left(\frac{x}{2^k}\right)\right) \ge \tau\left(e, \exp\left(\frac{x}{2^{k+1}}\right)\right) + \tau\left(\exp\left(\frac{x}{2^{k+1}}\right), \exp\left(\frac{x}{2^k}\right)\right) = 2\tau\left(e, \exp\left(\frac{x}{2^{k+1}}\right)\right) \ge 0,
\]
which implies $x \in C_k$. Let
\[
C_\infty := \bigcap_{k=1}^\infty C_k.
\]
The set $C_\infty$ is a nonempty compact set in $B_r(0)$ containing $0$.

We prove that $C_\infty \neq \{0\}$. Since $\tau(e, q) > 0$ and $q \neq e$, by property 6 there exists a continuous $\tau$-admissible curve $\gamma : [0, 1] \to G$ such that $\gamma(0) = e$ and $\gamma(1) = q$. Since $\gamma(0) = e \in B^G_{\frac{r}{2^k}}(e)$ and $\gamma(1) = q \notin B^G_r(e) \supset B^G_{\frac{r}{2^k}}(e)$, for each $k \ge 1$, by the connectedness of the interval $[0, 1]$, there exists a time $s_k \in (0, 1)$ such that $g_k := \gamma(s_k) \in \partial B^G_{\frac{r}{2^k}}(e)$. Due to the $\tau$-admissibility of the curve $\gamma$, we have $\tau(e, g_k) = \tau(\gamma(0), \gamma(s_k)) \ge 0$; therefore, $g_k \in M_{\frac{r}{2^k}}$, and the vector $a_k := 2^k \exp^{-1}(g_k)$ belongs to $C_k$ and satisfies $|a_k| = r$. By the compactness of the sphere $\partial B_r(0)$, one can extract a convergent subsequence $a_{k_j} \to a_*$ as $j \to \infty$ from the sequence $\{a_k\}_{k=1}^\infty$, and for its limit we have $|a_*| = r > 0$. Fix an arbitrary index $m \ge 1$. Then for all sufficiently large $j$, the inequality $k_j \ge m$ holds, and by the nestedness of the sets $C_k$, we have $a_{k_j} \in C_{k_j} \subset C_m$. Since the set $C_m$ is closed, taking the limit as $j \to \infty$ yields $a_* \in C_m$. Because $m \ge 1$ is arbitrary, we obtain $a_* \in \bigcap_{m=1}^\infty C_m = C_\infty$, which implies $C_\infty \neq \{0\}$.

We prove the convexity of the compact set $C_\infty$. Because $C_\infty$ is closed, it suffices to verify that for any $a, b \in C_\infty$, their half-sum $\frac{a+b}{2}$ belongs to $C_\infty$. Fix an arbitrary index $m \ge 1$. For all $k > m$, for the elements $g_k := \exp\left(\frac{a}{2^k}\right)$ and $h_k := \exp\left(\frac{b}{2^k}\right)$, we have $\tau(e, g_k) \ge 0$ and $\tau(g_k, g_k h_k) = \tau(e, h_k) \ge 0$, which yields $\tau(e, g_k h_k) \ge 0$. Write the product $g_k h_k$ in the form
\[
g_k h_k = \exp\left(\frac{a}{2^k} * \frac{b}{2^k}\right) = \exp\left(\frac{w_k}{2^{k-1}}\right), \qquad \text{where } w_k := 2^{k-1}\left(\frac{a}{2^k} * \frac{b}{2^k}\right) = \frac{a+b}{2} + 2^{k-1}z\left(\frac{a}{2^k}, \frac{b}{2^k}\right).
\]
From the estimate $|z(x, y)| \le |x| \cdot |y|$ it follows that $\left|w_k - \frac{a+b}{2}\right| \le \frac{r^2}{2^{k+1}} \to 0$ as $k \to \infty$, i.e., $w_k \to \frac{a+b}{2}$. At the same time, for all $k > m$, the reverse triangle inequality yields
\begin{equation}\label{eq:conv_ineq}
\tau\left(e, \exp\left(\frac{w_k}{2^m}\right)\right) \ge 2^{k-1-m}\tau\left(e, \exp\left(\frac{w_k}{2^{k-1}}\right)\right) = 2^{k-1-m}\tau(e, g_k h_k) \ge 0.
\end{equation}
Passing to the upper limit as $k \to \infty$ and using the upper semicontinuity of $\tau$ on $U \times U$ (property 4), we obtain $\tau\left(e, \exp\left(\frac{a+b}{2 \cdot 2^m}\right)\right) \ge 0$. Since $\left|\frac{a+b}{2}\right| \le r$, it follows that $\frac{a+b}{2} \in C_m$ for each $m \ge 1$, that is, $\frac{a+b}{2} \in C_\infty$. Thus, the set $C_\infty$ is convex.

We define the cone
\[
C_0^+ := \bigcup_{\mu > 0} \mu C_\infty.
\]

We prove the closedness of the cone $C_0^+$. Let a sequence $\{x_n\}_{n=1}^\infty \subset C_0^+$ converge to a vector $x \in L$. If $x = 0$, then $x \in C_\infty \subset C_0^+$. Suppose $x \neq 0$. Multiplying, if necessary, all terms of the sequence and the limit by a positive number $\frac{r}{2|x|}$, we may assume that $|x| = \frac{r}{2}$ and $|x_n| \le r$ for all sufficiently large $n$. For each such $n$, we have $x_n = \mu_n a_n$, where $a_n \in C_\infty$ and $\mu_n > 0$. Choose an integer $j_n \ge 1$ such that $\mu_n \le 2^{j_n}$. Then the scaled vector $\frac{x_n}{2^{j_n}} = \frac{\mu_n}{2^{j_n}}a_n$ belongs to the convex set $C_\infty$, which implies that for any $m \ge 1$, by the reverse triangle inequality,
\[
\tau\left(e, \exp\left(\frac{x_n}{2^m}\right)\right) \ge 2^{j_n}\tau\left(e, \exp\left(\frac{x_n}{2^{m+j_n}}\right)\right) \ge 0.
\]
Since $|x_n| \le r$, this means that $x_n \in C_m$ for all $m \ge 1$, i.e., $x_n \in C_\infty$ (in particular, $C_0^+ \cap B_r(0) = C_\infty$). By the closedness of the compact set $C_\infty$, passing to the limit as $n \to \infty$ yields $x \in C_\infty \subset C_0^+$. Thus, the cone $C_0^+$ is closed.

We prove the convexity of the cone $C_0^+$. Let $x, y \in C_0^+$; then by definition $x = \mu_1 a$ and $y = \mu_2 b$ for some $a, b \in C_\infty$ and $\mu_1, \mu_2 > 0$. Let $\mu := \max(\mu_1, \mu_2) > 0$. Since $0 \in C_\infty$ and the set $C_\infty$ is convex, the scaled vectors $\frac{x}{\mu} = \frac{\mu_1}{\mu}a$ and $\frac{y}{\mu} = \frac{\mu_2}{\mu}b$ belong to $C_\infty$, and therefore their half-sum $\frac{x+y}{2\mu}$ also lies in $C_\infty$. Scaling it back by a factor of $\mu$, we obtain $\frac{x+y}{2} = \mu\left(\frac{x+y}{2\mu}\right) \in C_0^+$. Consequently, the cone $C_0^+$ is convex.

We verify that the cone is pointed: $C_0^+ \cap (-C_0^+) = \{0\}$. If there existed a nonzero vector $w \in C_0^+ \cap (-C_0^+)$, then for the normalized vector $a := \frac{r}{|w|}w$, both vectors $a, -a$ would lie in $C_0^+ \cap B_r(0) = C_\infty \subset C_1$. Then for the element $g := \exp\left(\frac{a}{2}\right) \neq e$, the inequalities $\tau(e, g) \ge 0$ and $\tau(g, e) = \tau\left(e, \exp\left(-\frac{a}{2}\right)\right) \ge 0$ would hold simultaneously, which contradicts property 2.

Now we prove that the linear span $L_0 = C_0^+ - C_0^+$ generates the entire Lie algebra $L$ via the commutator operation. Let $\widetilde{L}_0 \subset L$ denote the smallest Lie subalgebra in $L$ containing $L_0$. We will need the following lemma.

\begin{lemma}\label{lem:3}
There exist a unit vector $\xi \in L$ ($|\xi| = 1$), a constant $c_0 \in \left(0, \frac{1}{2}\right)$, a radius $\rho := \frac{c_0 r}{4} \in \left(0, \frac{r}{8}\right)$, and an integer $k_1 \ge 1$ such that:
\begin{enumerate}
\item for all $u \in C_0^+$, the inequality $\langle \xi, u \rangle \ge 2c_0 |u|$ holds;
\item for all $k \ge k_1$ and for any vector $\tilde{x} \in A_{\frac{r}{2^k}}$ with norm $|\tilde{x}| = \frac{r}{2^k}$, the following estimate is valid:
\begin{equation}\label{eq:xi_step}
\langle \xi, \tilde{x} \rangle \ge c_0 |\tilde{x}| = c_0 \frac{r}{2^k};
\end{equation}
\item for any element $g \in B_\rho^G(e)$ that can be connected to the identity $e$ by a continuous $\tau$-admissible curve $\gamma : [0, 1] \to B_\rho^G(e)$ ($\gamma(0) = e$, $\gamma(1) = g$), the following representation holds for each $k \ge k_1$:
\[
g = h_1 h_2 \dots h_{N_k}, \qquad \text{where } h_i \in M_{\frac{r}{2^k}} \quad (N_k \le 2^k).
\]
\end{enumerate}
\end{lemma}

\begin{proof}[Proof of Lemma~\ref{lem:3}]
Since $C_0^+$ is a nontrivial, convex, closed, pointed cone ($C_0^+ \cap (-C_0^+) = \{0\}$), the intersection $S := C_0^+ \cap \partial B_1(0)$ is a nonempty compact set, and its convex hull $\operatorname{conv}(S)$ is compact and does not contain the vector $0$ (otherwise $0$ could be represented as a convex combination of nonzero vectors from $C_0^+$, which contradicts the pointedness of the cone). By the Hahn--Banach separation theorem applied to the convex compact set $\operatorname{conv}(S)$ and the point $\{0\}$, there exist a unit vector $\xi \in L$ ($|\xi| = 1$) and a number $c_0 \in \left(0, \frac{1}{2}\right)$ such that $\langle \xi, u \rangle \ge 2c_0$ for all $u \in S$. By positive homogeneity, this implies $\langle \xi, u \rangle \ge 2c_0 |u|$ for all $u \in C_0^+$.

For each $k \ge 1$, consider the closed subset of the compact sphere $\partial B_r(0)$:
\[
D_k := \left\{a \in C_k \cap \partial B_r(0) : \langle \xi, a \rangle \le c_0 r\right\}.
\]
Since $C_{k+1} \subset C_k$, the sequence of sets $\{D_k\}$ is nested: $D_{k+1} \subset D_k$. Let us prove by contradiction that there exists an integer $k_1 \ge 1$ for which $D_{k_1} = \varnothing$ (and then, by nestedness, $D_k = \varnothing$ for all $k \ge k_1$). Suppose the contrary: let $D_k \neq \varnothing$ for all $k \ge 1$. Then for each $k \ge 1$, we choose an element $a_k \in D_k$. Since the sphere $\partial B_r(0)$ is compact in the finite-dimensional space $L$, one can extract a convergent subsequence $a_{k_j} \to a_* \in \partial B_r(0)$ from the sequence $\{a_k\}$. Fix an arbitrary integer $m \ge 1$. Since $k_j \to \infty$, there is an index $j_0$ such that $k_j \ge m$ for all $j \ge j_0$. By the nestedness $D_{k_j} \subset D_m$, we have $a_{k_j} \in D_m$ for all $j \ge j_0$. Since the set $D_m$ is closed, passing to the limit as $j \to \infty$ gives $a_* \in D_m$. Because $m \ge 1$ was chosen arbitrarily, we obtain $a_* \in \bigcap_{m=1}^\infty D_m$. On the one hand, this means that $\langle \xi, a_* \rangle \le c_0 r$. On the other hand, $a_* \in \bigcap_{m=1}^\infty C_m = C_\infty \subset C_0^+$ and $|a_*| = r$, which implies, by the property of the vector $\xi$, that $\langle \xi, a_* \rangle \ge 2c_0 |a_*| = 2c_0 r > c_0 r$. The obtained contradiction proves the existence of an integer $k_1 \ge 1$ such that $D_k = \varnothing$ for all $k \ge k_1$.

This means that for any $a \in C_k$ with $|a| = r$, the inequality $\langle \xi, a \rangle > c_0 r$ holds for all $k \ge k_1$. Substituting $a = 2^k \tilde{x}$, where $\tilde{x} \in A_{\frac{r}{2^k}}$ and $|\tilde{x}| = \frac{r}{2^k}$ (so that $a \in C_k$ and $|a| = r$), and dividing both sides of the inequality $\langle \xi, 2^k \tilde{x} \rangle > c_0 r$ by $2^k$, we obtain the estimate~\eqref{eq:xi_step}.

Let $\rho := \frac{c_0 r}{4} \in (0, r/8)$. Note that for any points $a, b \in B_\rho^G(e)$, the vectors $z_a = \exp^{-1}(a)$ and $z_b = \exp^{-1}(b)$ have norms not exceeding $\rho$; therefore, by Lemma~\ref{lember:2}, the increment vector $\tilde{x} = (-z_a) * z_b = \exp^{-1}(a^{-1}b)$ is well-defined and satisfies the estimate $|\tilde{x}| \le 2\rho + \rho^2 < 3\rho < r < \lambda$.

Now let $\gamma : [0, 1] \to B_\rho^G(e)$ be a continuous $\tau$-admissible curve with $\gamma(0) = e$ and $\gamma(1) = g$, and let $k \ge k_1$. We inductively construct time instances $0 = s_0 < s_1 < \dots < s_{N_k} = 1$:
\begin{itemize}
\item if the time $s_{i-1} < 1$ has already been constructed and for all $s \in [s_{i-1}, 1]$ the inclusion $\gamma(s_{i-1})^{-1}\gamma(s) \in B^G_{\frac{r}{2^k}}(e)$ holds, we set $N_k := i$ and $s_{N_k} := 1$;
\item otherwise, we set
\[
s_i := \min \left\{s \in [s_{i-1}, 1] : \gamma(s_{i-1})^{-1}\gamma(s) \in \partial B^G_{\frac{r}{2^k}}(e)\right\}.
\]
\end{itemize}
Due to the uniform continuity of $\gamma$, the process will terminate in a finite number of steps $N_k$. Moreover, the elements $h_i := \gamma(s_{i-1})^{-1}\gamma(s_i)$ satisfy the conditions: $\tau(e, h_i) = \tau(\gamma(s_{i-1}), \gamma(s_i)) \ge 0$ and $h_i \in B^G_{\frac{r}{2^k}}(e)$, which means $h_i \in M_{\frac{r}{2^k}}$ for all $i = 1, \dots, N_k$. For the first $N_k - 1$ steps, the vectors $\tilde{x}_i := \exp^{-1}(h_i) \in A_{\frac{r}{2^k}}$ have norm $|\tilde{x}_i| = \frac{r}{2^k}$.

Let us estimate the number $N_k$ from above. Denote $z_i := \exp^{-1}(\gamma(s_i)) \in B_\rho(0)$. Then $z_i = z_{i-1} * \tilde{x}_i = z_{i-1} + \tilde{x}_i + z(z_{i-1}, \tilde{x}_i)$, where by Lemma~\ref{lember:1} we have $|z(z_{i-1}, \tilde{x}_i)| \le |z_{i-1}| \cdot |\tilde{x}_i| \le \rho |\tilde{x}_i|$. For each complete step $i \in \{1, \dots, N_k - 1\}$, taking into account~\eqref{eq:xi_step} and the inequality $\rho = \frac{c_0 r}{4} < \frac{c_0}{2}$, we have:
\[
\langle \xi, z_i \rangle - \langle \xi, z_{i-1} \rangle = \langle \xi, \tilde{x}_i \rangle + \big\langle \xi, z(z_{i-1}, \tilde{x}_i) \big\rangle \ge c_0|\tilde{x}_i| - \rho|\tilde{x}_i| \ge \frac{c_0}{2}|\tilde{x}_i| = \frac{c_0 r}{2 \cdot 2^k}.
\]
Summing these inequalities over $i = 1, \dots, N_k - 1$ and noting that $z_0 = 0$ and $|z_{N_k - 1}| \le \rho$, we obtain:
\[
(N_k - 1)\frac{c_0 r}{2 \cdot 2^k} \le \langle \xi, z_{N_k - 1} \rangle - \langle \xi, z_0 \rangle \le |z_{N_k - 1}| \le \rho = \frac{c_0 r}{4}.
\]
Dividing both sides by $\frac{c_0 r}{2 \cdot 2^k} > 0$, we find $N_k - 1 \le 2^{k-1}$, which yields $N_k \le 2^{k-1} + 1 \le 2^k$.
\end{proof}

We will need one more lemma.

\begin{lemma}\label{lem:4}
For the radius $\rho > 0$ from Lemma~\ref{lem:3}, there exists a point $q' \in B_\rho^G(e)$ with $\tau(e, q') > 0$ such that the set $W_G := \operatorname{int} I(e, q')$ is nonempty, open in $G$, and its closure is contained in the closed chronological diamond $J(e, q') \subset B_\rho^G(e)$.
\end{lemma}

\begin{proof}[Proof of Lemma~\ref{lem:4}]
1. Applying Lemma~\ref{lem:1} to the open neighborhood of the identity $\operatorname{int} B_\rho^G(e)$, we find points $p_1, q_1 \in \operatorname{int} B_\rho^G(e)$ such that $e \in \operatorname{int} I(p_1, q_1) \subset I(p_1, q_1) \subset B_\rho^G(e)$, and for any $p, q' \in I(p_1, q_1)$ the closed chronological diamond $J(p, q')$ is compact and lies entirely in $I(p_1, q_1) \subset B_\rho^G(e)$.

2. Note that for any point $q \in G$, the condition $e \in \operatorname{int} I^-(q)$ implies $q \in \operatorname{int} I^+(e)$. Indeed, if $e \in \operatorname{int} I^-(q)$, then there exists an open neighborhood $V_e \subset G$ of the identity $e$ such that $V_e \subset I^-(q)$, i.e., $\tau(y, q) > 0$ for all $y \in V_e$. Since the mapping $y \mapsto y^{-1}q$ is a homeomorphism of the group $G$ onto itself mapping $e$ to $q$, the set $U_q := \{y^{-1}q : y \in V_e\}$ is an open neighborhood of the point $q$. Moreover, for each point $z = y^{-1}q \in U_q$ ($y \in V_e$), by the left-invariance of $\tau$ (property 1) we have
\[
\tau(e, z) = \tau(e, y^{-1}q) = \tau(ye, y(y^{-1}q)) = \tau(y, q) > 0,
\]
which means $U_q \subset I^+(e)$, whence $q \in \operatorname{int} I^+(e)$.

3. Applying property 5 to the point $e$ and its open neighborhood $\operatorname{int} I(p_1, q_1)$, we find points $p', q' \in \operatorname{int} I(p_1, q_1)$ such that
\[
e \in \operatorname{int} I(p', q') \subset I(p', q') \subset \operatorname{int} I(p_1, q_1).
\]
In particular, $q' \in I(p_1, q_1) \subset B_\rho^G(e)$ and $\tau(e, q') > 0$. 

Applying property 5 once again to the point $e$ and its open neighborhood $\operatorname{int} I(p', q')$, we find points $p'', q'' \in \operatorname{int} I(p', q')$ for which
\[
e \in \operatorname{int} I(p'', q'') \subset I(p'', q'') \subset \operatorname{int} I(p', q').
\]
Since $e \in \operatorname{int} I(p'', q'') = \operatorname{int} I^+(p'') \cap \operatorname{int} I^-(q'') \subset \operatorname{int} I^-(q'')$, it follows from step 2 that $q'' \in \operatorname{int} I^+(e)$. On the other hand, by the choice of the point $q''$, we have $q'' \in \operatorname{int} I(p', q') \subset \operatorname{int} I^-(q')$. Consequently,
\[
q'' \in \operatorname{int} I^+(e) \cap \operatorname{int} I^-(q') = \operatorname{int} I(e, q').
\]
Thus, the open set $W_G := \operatorname{int} I(e, q')$ contains the point $q''$ and is therefore nonempty. Finally, since $e, q' \in I(p_1, q_1)$, by step 1 the closed chronological diamond $J(e, q')$ is compact (and hence closed in $G$), which yields $\overline{W}_G \subset J(e, q') \subset B_\rho^G(e)$.
\end{proof}

By Lemma~\ref{lem:4}, the set $W_G = \operatorname{int} I(e, q') \subset J(e, q') \subset B_\rho^G(e)$ is nonempty and open in $G$; therefore, its preimage $W := \exp^{-1}(W_G) \subset B_\rho(0)$ is a nonempty open subset in the Lie algebra $L$. It suffices to prove that $W \subset \widetilde{L}_0$: since the Lie subalgebra $\widetilde{L}_0$ is a vector subspace in $L$ containing the nonempty open set $W$ in $L$, this implies $\widetilde{L}_0 = L$.

Let $x \in W$. Then $g := \exp(x) \in W_G \subset I(e, q')$; in particular, $\tau(e, g) > 0$. By property 6, there exists a continuous $\tau$-admissible curve $\gamma : [0, 1] \to G$ connecting $e$ to $g$. For any $s \in [0, 1]$, we have $\tau(e, \gamma(s)) \ge 0$ and $\tau(\gamma(s), q') \ge \tau(\gamma(s), g) + \tau(g, q') > 0$, so the entire image $\gamma([0, 1])$ lies in $J(e, q') \subset B_\rho^G(e)$.

Take an arbitrary $\varepsilon \in \left(0, \min\left(1, \frac{\lambda}{14 r e^r}\right)\right)$. Since $C_\infty = \bigcap_{k=1}^\infty C_k$ is the intersection of a decreasing sequence of compact sets in $B_r(0)$, for all sufficiently large $k \ge 1$ the inclusion $C_k \subset C_\infty + U(0, \varepsilon r)$ holds. Dividing by $2^k$, we obtain:
\begin{equation}\label{eq:A_incl}
A_{\frac{r}{2^k}} = \frac{1}{2^k}C_k \subset \frac{1}{2^k}\big(C_\infty + U(0, \varepsilon r)\big).
\end{equation}

By Lemma~\ref{lem:3}, the element $g = \exp(x)$ can be represented in the form $g = h_1 h_2 \dots h_{N_k}$, where $h_i \in M_{\frac{r}{2^k}}$ and $N_k \le 2^k$. Setting $\tilde{x}_i := \exp^{-1}(h_i) \in A_{\frac{r}{2^k}}$, in exponential coordinates we have:
\[
x = \tilde{x}_1 * \tilde{x}_2 * \dots * \tilde{x}_{N_k}.
\]
Due to the inclusion~\eqref{eq:A_incl} and the inequality $N_k \le 2^k$, each vector $\tilde{x}_i$ can be represented as
\[
\tilde{x}_i = x_i + y_i, \qquad \text{where } x_i \in \frac{1}{2^k}C_\infty \subset B_{\frac{r}{2^k}}(0), \quad |x_i| \le \frac{r}{2^k} \le \frac{r}{N_k}, \quad |y_i| < \frac{\varepsilon r}{2^k} \le \frac{\varepsilon r}{N_k}.
\]

Applying Lemma~\ref{lember:6} with $\eta = r$ and $j = N_k$ (note that the restriction on $\varepsilon$ from Lemma~\ref{lember:6} is satisfied by our choice), we obtain:
\begin{equation}\label{eq:x_decomp}
x = (x_1 + y_1) * \dots * (x_{N_k} + y_{N_k}) = (x_1 * x_2 * \dots * x_{N_k}) * y, \qquad \text{where } |y| < 14\varepsilon r e^r.
\end{equation}
Since $x_i \in \frac{1}{2^k}C_\infty \subset C_0^+ \subset L_0 \subset \widetilde{L}_0$, and the Lie subalgebra $\widetilde{L}_0$ is closed with respect to the operation $*$ (due to the Baker--Campbell--Hausdorff formula), the vector $w_k := x_1 * x_2 * \dots * x_{N_k}$ belongs to the Lie subalgebra $\widetilde{L}_0$, and $|w_k| < e^r - 1 < 1$ by Lemma~\ref{lember:5}. By Lemma~\ref{lember:1}, estimate~\eqref{eq:x_decomp} implies:
\[
|x - w_k| = |w_k * y - w_k| \le |y| + |w_k| \cdot |y| < 2 \cdot 14\varepsilon r e^r = 28\varepsilon r e^r.
\]

Since $\varepsilon > 0$ can be chosen arbitrarily small, and the finite-dimensional subspace $\widetilde{L}_0 \subset L$ is closed, as $\varepsilon \to 0$ we obtain $x \in \widetilde{L}_0$. Thus, $W \subset \widetilde{L}_0$, which yields $\widetilde{L}_0 = L$.

Before proceeding to the definition of the anti-norm, we prove one more lemma.

\begin{lemma}\label{lem:5}
The following statements hold:
\begin{enumerate}
\item A vector $u \in L$ belongs to the cone $C_0^+$ if and only if there exists $\delta > 0$ such that $\tau(e, \exp(tu)) \ge 0$ for all $t \in [0, \delta]$. Moreover, for any $u \in C_0^+$, the inequality $\tau(e, \exp(tu)) \ge 0$ holds for all $t \ge 0$.
\item For each vector $u \in C_0^+$, the function $f_u(t) := \tau(e, \exp(tu))$ on the half-line $[0, +\infty)$ is non-negative, satisfies the reverse triangle inequality ($f_u(t + s) \ge f_u(t) + f_u(s)$ for all $t, s \ge 0$), and is non-decreasing. Furthermore, there exists a finite limit
\begin{equation}\label{eq:nu0_def}
\nu_0(u) := \lim_{k \to \infty} 2^k \tau\!\left(e, \exp\!\left(\frac{u}{2^k}\right)\right) = \lim_{t \to 0^+} \frac{\tau(e, \exp(tu))}{t} = \inf_{t > 0} \frac{\tau(e, \exp(tu))}{t},
\end{equation}
satisfying the bound
\begin{equation}\label{eq:nu0_bound}
0 \le \nu_0(u) \le \frac{\tau(e, q)}{r}|u|.
\end{equation}
\end{enumerate}
\end{lemma}
\begin{proof}[Proof of Lemma \ref{lem:5}]
1. For $u = 0$, we have $0 \in C_0^+$ and $\tau(e, \exp(0)) = \tau(e, e) = 0 \ge 0$ for all $t \ge 0$. Let $u \in C_0^+ \setminus \{0\}$ and $t > 0$. Set $u_0 := \frac{r}{|u|}u$. Since $|u_0| = r$ and $C_0^+ \cap B_r(0) = C_\infty$, for any $\alpha \in (0, 1]$ the vector $\alpha u_0$ belongs to $C_\infty \subset C_1$. Choose an integer $k \ge 1$ large enough so that the number $\alpha := \frac{t|u|}{2^{k-1}r}$ belongs to the half-open interval $(0, 1]$. Then $\alpha u_0 \in C_1$, whence by the definition of the set $C_1$ we obtain
\[
\tau\!\left(e, \exp\!\left(\frac{tu}{2^k}\right)\right) = \tau\!\left(e, \exp\!\left(\frac{\alpha u_0}{2}\right)\right) \ge 0.
\]
Since $\exp(tu) = \left(\exp\!\left(\frac{tu}{2^k}\right)\right)^{2^k}$, by the reverse triangle inequality and the left-invariance of $\tau$, we have
\[
\tau(e, \exp(tu)) \ge 2^k \tau\!\left(e, \exp\!\left(\frac{tu}{2^k}\right)\right) \ge 0.
\]
Conversely, let $u \in L \setminus \{0\}$ and suppose $\tau(e, \exp(tu)) \ge 0$ for all $t \in [0, \delta]$, where $\delta > 0$. Set $u_0 := \frac{r}{|u|}u \in \partial B_r(0)$. For all sufficiently large integers $k \ge 1$ such that $t_k := \frac{r}{2^k |u|} \le \delta$, we have
\[
\tau\!\left(e, \exp\!\left(\frac{u_0}{2^k}\right)\right) = \tau(e, \exp(t_k u)) \ge 0,
\]
which means $u_0 \in C_k$. By the nesting $C_{k+1} \subset C_k$, it follows that $u_0 \in C_k$ for all $k \ge 1$. Hence, $u_0 \in \bigcap_{k=1}^\infty C_k = C_\infty \subset C_0^+$, and therefore $u = \frac{|u|}{r}u_0 \in C_0^+$.
        
2. Let $u \in C_0^+$. By part 1, the function $f_u(t) = \tau(e, \exp(tu))$ is non-negative for all $t \ge 0$. For any $t, s \ge 0$, from the identity $\exp((t + s)u) = \exp(tu)\exp(su)$, the left-invariance of $\tau$, and the reverse triangle inequality for $\tau$, we obtain the reverse triangle inequality for the function $f_u$:
\[
f_u(t + s) = \tau(e, \exp(tu)\exp(su)) \ge \tau(e, \exp(tu)) + \tau(\exp(tu), \exp(tu)\exp(su)) = f_u(t) + f_u(s).
\]
Since $f_u(s) \ge 0$, this implies $f_u(t + s) \ge f_u(t)$, meaning the function $f_u$ is non-decreasing on $[0, +\infty)$.
For $u = 0$, we have $f_0(t) \equiv 0$, and relations \eqref{eq:nu0_def} and \eqref{eq:nu0_bound} are trivially satisfied. Let $u \in C_0^+ \setminus \{0\}$. For arbitrary $t > 0$ and $s \in (0, t)$, set $m := \lfloor t/s \rfloor \ge 1$ and $r_s := t - ms \in [0, s)$. By the reverse triangle inequality and the non-negativity of $f_u$, we have
\[
f_u(t) = f_u(ms + r_s) \ge m f_u(s) + f_u(r_s) \ge m f_u(s),
\]
whence dividing by $ms > 0$ yields
\[
\frac{f_u(s)}{s} \le \frac{t}{ms} \cdot \frac{f_u(t)}{t} = \frac{t/s}{\lfloor t/s \rfloor} \cdot \frac{f_u(t)}{t}.
\]
Taking the limit superior as $s \to 0^+$ and noting that $\lim_{s \to 0^+} \frac{t/s}{\lfloor t/s \rfloor} = 1$, we find
\[
\varlimsup_{s \to 0^+} \frac{f_u(s)}{s} \le \frac{f_u(t)}{t} \qquad \forall t > 0.
\]
Taking the infimum over all $t > 0$ on the right-hand side yields the inequality
\[
\varlimsup_{s \to 0^+} \frac{f_u(s)}{s} \le \inf_{t > 0} \frac{f_u(t)}{t}.
\] 
On the other hand, since for any $\delta > 0$ it holds that
\[
\inf_{t > 0} \frac{f_u(t)}{t} \le \inf_{s \in (0, \delta)} \frac{f_u(s)}{s},
\]
taking the supremum over $\delta > 0$ gives
\[
\inf_{t > 0} \frac{f_u(t)}{t} \le \sup_{\delta > 0} \inf_{s \in (0, \delta)} \frac{f_u(s)}{s} = \varliminf_{s \to 0^+} \frac{f_u(s)}{s}.
\] 
Thus, we obtain the chain of inequalities
\[
\varlimsup_{s \to 0^+} \frac{f_u(s)}{s} \le \inf_{t > 0} \frac{f_u(t)}{t} \le \varliminf_{s \to 0^+} \frac{f_u(s)}{s},
\]
which proves the existence of the limit $\lim_{t \to 0^+} \frac{f_u(t)}{t} = \inf_{t > 0} \frac{f_u(t)}{t}$ and, in particular, equality \eqref{eq:nu0_def} upon substituting $t = 2^{-k}$. Finally, setting $t_0 := \frac{r}{|u|} > 0$, we observe that $|t_0 u| = r$, and therefore $\exp(t_0 u) \in B_r^G(e) \subset I(p, q)$. By the reverse triangle inequality, $\tau(e, q) \ge \tau(e, \exp(t_0 u)) + \tau(\exp(t_0 u), q) > f_u(t_0)$, whence
\[
0 \le \nu_0(u) = \inf_{t > 0} \frac{f_u(t)}{t} \le \frac{f_u(t_0)}{t_0} \le \frac{\tau(e, q)}{r}|u|.
\]
\end{proof}

We define the anti-norm $\nu_0 : L \to [0, +\infty) \cup \{-\infty\}$ by the formula
\[
\nu_0(u) := \begin{cases}
\lim\limits_{t \to 0^+} \dfrac{\tau(e, \exp(tu))}{t} = \inf\limits_{t > 0} \dfrac{\tau(e, \exp(tu))}{t}, & \text{if } u \in C_0^+, \\[2ex]
-\infty, & \text{if } u \in L \setminus C_0^+.
\end{cases}
\]

Let us verify that the function $\nu_0$ satisfies all the properties of an anti-norm from the definition of a sub-Lorentzian structure.

We prove positive homogeneity. For any $u \in C_0^+$ and any $\mu > 0$, we have $\mu u \in C_0^+$; thus, making the substitution $s = \mu t$ in formula~\eqref{eq:nu0_def}, we obtain:
\[
\nu_0(\mu u) = \lim_{t \to 0^+} \frac{\tau(e, \exp(t\mu u))}{t} = \mu \lim_{s \to 0^+} \frac{\tau(e, \exp(su))}{s} = \mu \nu_0(u).
\]
For $\mu = 0$ and $u \in C_0^+$, the equality $\nu_0(0 \cdot u) = \nu_0(0) = 0 = 0 \cdot \nu_0(u)$ trivially holds, and for $u \in L \setminus C_0^+$ and $\mu > 0$, we have $\mu u \in L \setminus C_0^+$ and $\nu_0(\mu u) = -\infty = \mu \nu_0(u)$.

We show that for any $u, v \in L$, the following inequality holds:
\[
\nu_0(u + v) \ge \nu_0(u) + \nu_0(v).
\]
If at least one of the vectors $u, v$ does not belong to $C_0^+$, then the right-hand side is $-\infty$, and the inequality is true; if $u = 0$ or $v = 0$, it is also true. Let $u, v \in C_0^+ \setminus \{0\}$. Due to the positive homogeneity of $\nu_0$, multiplying $u$ and $v$ by a suitable common positive factor, we may assume without loss of generality that $u, v \in C_0^+ \cap B_r(0) = C_\infty$. Then for each $k \ge 1$, the vectors $\frac{u}{2^k}$ and $\frac{v}{2^k}$ lie in $A_{\frac{r}{2^k}} \subset B_{\frac{r}{2^k}}(0)$, so by Lemma~\ref{lember:2}, the element
\[
\tilde{w}_k := \left(\frac{u}{2^k}\right) * \left(\frac{v}{2^k}\right) = \exp^{-1}\!\left(\exp\!\left(\frac{u}{2^k}\right)\exp\!\left(\frac{v}{2^k}\right)\right) \in V
\]
is well-defined. Set $w_k := 2^{k-1}\tilde{w}_k$. By Lemma~\ref{lember:1}, we have
\[
\left|w_k - \frac{u + v}{2}\right| = 2^{k-1}\left|z\!\left(\frac{u}{2^k}, \frac{v}{2^k}\right)\right| \le 2^{k-1}\frac{|u|}{2^k}\frac{|v|}{2^k} = \frac{|u| \cdot |v|}{2^{k+1}} \le \frac{r^2}{2^{k+1}} \to 0 \qquad (k \to \infty),
\]
meaning $w_k \to \frac{u + v}{2}$ as $k \to \infty$ (where $|w_k| \le r + \frac{r^2}{4} < 2r < \lambda$). Fix an arbitrary integer $m \ge 1$. For all $k \ge m + 1$, successively applying the reverse triangle inequality and the left-invariance of $\tau$, we obtain:
\begin{align*}
2^m \tau\!\left(e, \exp\!\left(\frac{w_k}{2^m}\right)\right) &\ge 2^m \cdot 2^{k-1-m} \tau\!\left(e, \exp\!\left(\frac{w_k}{2^{k-1}}\right)\right) = 2^{k-1}\tau(e, \exp(\tilde{w}_k)) = \\
&= 2^{k-1}\tau\!\left(e, \exp\!\left(\frac{u}{2^k}\right)\exp\!\left(\frac{v}{2^k}\right)\right) \ge \\
&\ge \frac{1}{2}\left(2^k \tau\!\left(e, \exp\!\left(\frac{u}{2^k}\right)\right) + 2^k \tau\!\left(e, \exp\!\left(\frac{v}{2^k}\right)\right)\right).
\end{align*}
Taking the limit superior as $k \to \infty$ in this inequality, utilizing the upper semicontinuity of $\tau$ on $U_0 \times U_0$ (Property 4) on the left-hand side and relation~\eqref{eq:nu0_def} on the right-hand side, we find:
\[
2^m \tau\!\left(e, \exp\!\left(\frac{u + v}{2^{m+1}}\right)\right) \ge \varlimsup_{k \to \infty} 2^m \tau\!\left(e, \exp\!\left(\frac{w_k}{2^m}\right)\right) \ge \frac{\nu_0(u) + \nu_0(v)}{2}.
\]
Letting $m \to \infty$, by formula~\eqref{eq:nu0_def} for the vector $\frac{u + v}{2} \in C_0^+$, we obtain
\[
\nu_0\!\left(\frac{u + v}{2}\right) \ge \frac{\nu_0(u) + \nu_0(v)}{2},
\]
whence, after multiplying by $2$, the positive homogeneity yields $\nu_0(u + v) \ge \nu_0(u) + \nu_0(v)$. Combined with positive homogeneity, this means that the function $\nu_0$ is concave on the convex cone $C_0^+$.

We verify upper semicontinuity on $L$ and continuity on $\operatorname{int}_{L_0}(C_0^+)$. Fix an arbitrary $R > 0$ and set $\delta := \lambda / R$. For each $t \in (0, \delta)$ and any $u \in U(0, R)$, we have $|tu| < \lambda$, and thus $\exp(tu) \in U \subset U_0$. Since the function $\tau(e, \cdot)$ is upper semicontinuous on $U_0$, for each fixed $t \in (0, \delta)$, the mapping $u \mapsto \frac{\tau(e, \exp(tu))}{t}$ is upper semicontinuous on the open ball $U(0, R)$. Due to equality~\eqref{eq:nu0_def}, for all $u \in C_0^+ \cap U(0, R)$, the following representation holds:
\[
\nu_0(u) = \inf_{t \in (0, \delta)} \frac{\tau(e, \exp(tu))}{t}.
\]
Since the infimum of any family of upper semicontinuous functions is upper semicontinuous, and $R > 0$ is arbitrary, the function $\nu_0$ is upper semicontinuous on the cone $C_0^+$. Being extended by the value $-\infty$ outside the closed cone $C_0^+$, it is also upper semicontinuous on the entire Lie algebra $L$. Furthermore, since any finite concave function defined on a convex set in a finite-dimensional space is continuous on its relative interior, the restriction $\nu_0|_{\operatorname{int}_{L_0}(C_0^+)}$ is continuous.

We prove strict positivity on $\operatorname{int}_{L_0}(C_0^+)$. We first show that there exists a vector $u_* \in C_\infty \subset C_0^+$ for which $\nu_0(u_*) > 0$. By Lemma~\ref{lem:4}, there is a point $q' \in B_\rho^G(e)$ with $\tau(e, q') > 0$ and $J(e, q') \subset B_\rho^G(e)$. By Property 6, there exists a continuous $\tau$-admissible curve $\gamma : [0, 1] \to G$ with $\gamma(0) = e$, $\gamma(1) = q'$, and a positive $\tau$-length $c := L_\tau(\gamma) > 0$. Since $\gamma$ is $\tau$-admissible, its image is contained in $J(e, q') \subset B_\rho^G(e)$. According to the construction from the proof of part 3 of Lemma~\ref{lem:3}, for each $k \ge k_1$ there exists a partition $0 = s_0 < s_1 < \dots < s_{N_k} = 1$ with $N_k \le 2^k$ steps, for which the elements $h_i := \gamma(s_{i-1})^{-1}\gamma(s_i)$ belong to $M_{\frac{r}{2^k}}$ ($i = 1, \dots, N_k$). By the definition of the $\tau$-length $L_\tau(\gamma)$ and the left-invariance of $\tau$, we have:
\[
c = L_\tau(\gamma) \le \sum_{i=1}^{N_k} \tau(\gamma(s_{i-1}), \gamma(s_i)) = \sum_{i=1}^{N_k} \tau(e, h_i).
\]
Since there are $N_k \le 2^k$ terms, by the pigeonhole principle, there exists an index $i_k \in \{1, \dots, N_k\}$ such that
\[
\tau(e, h_{i_k}) \ge \frac{c}{N_k} \ge \frac{c}{2^k}.
\]
Set $u_k := 2^k \exp^{-1}(h_{i_k})$. Since $h_{i_k} \in M_{\frac{r}{2^k}}$, the vector $u_k$ belongs to the compact set $C_k \subset B_r(0)$ and satisfies the inequality
\[
2^k \tau\!\left(e, \exp\!\left(\frac{u_k}{2^k}\right)\right) = 2^k \tau(e, h_{i_k}) \ge c > 0.
\]
By the compactness of the ball $B_r(0)$, we extract from the sequence $\{u_k\}_{k \ge k_1}$ a convergent subsequence $u_{k_j} \to u_* \in B_r(0)$ as $j \to \infty$. Just as before, the nesting of the closed sets $C_k$ implies that $u_* \in \bigcap_{m=1}^\infty C_m = C_\infty \subset C_0^+$. Fix an arbitrary integer $m \ge 1$. For all $k_j \ge m$, the reverse triangle inequality yields
\[
2^m \tau\!\left(e, \exp\!\left(\frac{u_{k_j}}{2^m}\right)\right) \ge 2^{k_j} \tau\!\left(e, \exp\!\left(\frac{u_{k_j}}{2^{k_j}}\right)\right) \ge c.
\]
Taking the limit superior as $j \to \infty$ and using the upper semicontinuity of $\tau$ on $U_0 \times U_0$, we obtain $2^m \tau\!\left(e, \exp\!\left(\frac{u_*}{2^m}\right)\right) \ge c$ for each $m \ge 1$. Letting $m \to \infty$, we find $\nu_0(u_*) \ge c > 0$.

Now let $u \in \operatorname{int}_{L_0}(C_0^+)$ be an arbitrary interior vector of the cone $C_0^+$ in the subspace $L_0 = C_0^+ - C_0^+$. Since $u_* \in C_0^+ \subset L_0$, by the definition of the relative interior, there exists a sufficiently small number $\varepsilon > 0$ such that the vector $v := u - \varepsilon u_*$ belongs to $C_0^+$. Then $u = v + \varepsilon u_*$, and by the concavity, positive homogeneity, and non-negativity of $\nu_0$ on $C_0^+$, we obtain:
\[
\nu_0(u) = \nu_0(v + \varepsilon u_*) \ge \nu_0(v) + \varepsilon \nu_0(u_*) \ge \varepsilon c > 0.
\]
Thus, the triple $(G, C_0^+, \nu_0)$ is a left-invariant sub-Lorentzian structure on the Lie group $G$.

Now, using this anti-norm $\nu_0$ and the cone $C_0^+$, we construct a sub-Lorentzian chronometer $\tau_c : G \times G \to [0, +\infty] \cup \{-\infty\}$ and show that $\tau(g, h) = \tau_c(g, h)$ for all $g, h \in G$. We prove the following lemma.

\begin{lemma}\label{lem:6}
There exists an integer $k_2 \ge k_1$ such that for any vector $x \in B_{\frac{r}{2^{k_2}}}(0)$ satisfying the condition $\tau(e, \exp(x)) \ge 0$, the following inequalities hold:
\begin{equation}\label{eq:local_cone_bound}
\langle \xi, x \rangle \ge c_0 |x|, \qquad 0 \le \tau(e, \exp(x)) \le \frac{2\tau(e, q)}{r}|x|,
\end{equation}
where $\xi$ and $c_0$ are taken from Lemma \ref{lem:3}.
\end{lemma}
\begin{proof}[Proof of Lemma \ref{lem:6}]
For $x = 0$, both inequalities \eqref{eq:local_cone_bound} turn into the trivial equality $0 = 0$. For each $k \ge 1$, consider the closed subset of the compact spherical shell $\left\{a \in L : \frac{r}{2} \le |a| \le r\right\}$:
\[
\widetilde{D}_k := \left\{a \in C_k : \frac{r}{2} \le |a| \le r,\ \langle \xi, a \rangle \le c_0 |a|\right\}.
\]
Since $C_{k+1} \subset C_k$, the sequence of sets $\{\widetilde{D}_k\}$ is nested: $\widetilde{D}_{k+1} \subset \widetilde{D}_k$. We prove by contradiction that there exists an integer $k_2 \ge k_1$ for which $\widetilde{D}_k = \varnothing$ for all $k \ge k_2$. If $\widetilde{D}_k \neq \varnothing$ for all $k \ge 1$, then, by choosing an element $a_k \in \widetilde{D}_k$ and extracting a convergent subsequence $a_{k_j} \to a_*$, the closedness and nestedness of the sets $\widetilde{D}_m$ imply $a_* \in \bigcap_{m=1}^\infty \widetilde{D}_m$. Then $\frac{r}{2} \le |a_*| \le r$, $\langle \xi, a_* \rangle \le c_0 |a_*|$, and simultaneously $a_* \in \bigcap_{m=1}^\infty C_m = C_\infty \subset C_0^+$. However, by part 1 of Lemma \ref{lem:3}, the inequality $\langle \xi, a_* \rangle \ge 2c_0 |a_*| > c_0 |a_*|$ holds for the vector $a_* \in C_0^+ \setminus \{0\}$---a contradiction. Consequently, there exists an integer $k_2 \ge k_1$ such that $\widetilde{D}_k = \varnothing$ for all $k \ge k_2$.

Now let $x \in B_{\frac{r}{2^{k_2}}}(0) \setminus \{0\}$ and $\tau(e, \exp(x)) \ge 0$. Choose the unique integer $m \ge k_2$ such that the double inequality
\[
\frac{r}{2^{m+1}} < |x| \le \frac{r}{2^m}
\]
holds, and set $a := 2^m x$. Then $\frac{r}{2} < |a| \le r$ and $\tau\!\left(e, \exp\!\left(\frac{a}{2^m}\right)\right) = \tau(e, \exp(x)) \ge 0$, meaning $a \in C_m$. Since $m \ge k_2$ and $\widetilde{D}_m = \varnothing$, we have $\langle \xi, a \rangle > c_0 |a|$. Dividing both sides of this inequality by $2^m$, we obtain the first inequality in \eqref{eq:local_cone_bound}: $\langle \xi, x \rangle \ge c_0 |x|$.

To prove the second inequality in \eqref{eq:local_cone_bound}, observe that since $\exp(a) = (\exp(x))^{2^m}$ and $\tau(e, \exp(x)) \ge 0$, the reverse triangle inequality and the left-invariance of $\tau$ yield
\[
0 \le 2^m \tau(e, \exp(x)) \le \tau(e, \exp(a)).
\]
At the same time, $|a| \le r$, so the point $\exp(a)$ belongs to the ball $B_r^G(e) \subset I(p, q)$, where $\tau(\exp(a), q) > 0$. Consequently, by the reverse triangle inequality, $\tau(e, \exp(a)) \le \tau(e, q) - \tau(\exp(a), q) < \tau(e, q)$. Given that $|x| > \frac{r}{2^{m+1}}$ implies $\frac{1}{2^m} < \frac{2|x|}{r}$, we finally obtain:
\[
0 \le \tau(e, \exp(x)) \le \frac{1}{2^m}\tau(e, \exp(a)) \le \frac{\tau(e, q)}{2^m} \le \frac{2\tau(e, q)}{r}|x|.
\]
\end{proof}

Let us first prove that $\tau(g, h) \le \tau_c(g, h)$ for all $g, h \in G$. If $g = h$, then by Property~2 we have $\tau(g, g) = 0$, and the constant curve $\gamma_0(t) \equiv g$ ($t \in [0, 1]$) is $C^+$-admissible and has sub-Lorentzian length $\ell(\gamma_0) = \int_0^1 \nu_0(0)\,dt = 0$, whence $\tau_c(g, g) \ge \ell(\gamma_0) = 0 = \tau(g, g)$.

Now let $g \neq h$. If $\tau(g, h) = -\infty$, then the inequality $\tau(g, h) \le \tau_c(g, h)$ trivially holds. Suppose $\tau(g, h) \ge 0$. By Property~6, the set $\Omega_\tau(g, h)$ is non-empty and $\tau(g, h) = \sup_{\gamma \in \Omega_\tau(g, h)} L_\tau(\gamma)$. Fix an arbitrary continuous $\tau$-admissible curve $\gamma \in \Omega_\tau(g, h)$ (parameterized on the interval $[0, 1]$ so that $\gamma(0) = g$, $\gamma(1) = h$) and we will show that $L_\tau(\gamma) \le \tau_c(g, h)$. Set
\[
\rho_0 := \min\left\{\rho, \frac{r}{4 \cdot 2^{k_2}}\right\} = \min\left\{\frac{c_0 r}{4}, \frac{r}{4 \cdot 2^{k_2}}\right\} > 0.
\]
By the uniform continuity of the curve $\gamma$ on the interval $[0, 1]$, we can choose a partition $0 = a_0 < a_1 < \dots < a_N = 1$ fine enough so that for each $j = 1, \dots, N$ and all $s \in [a_{j-1}, a_j]$, the inclusion $\gamma(a_{j-1})^{-1}\gamma(s) \in B_{\rho_0}^G(e)$ holds (moreover, if on some segment $\gamma(a_{j-1}) = \gamma(a_j)$, then by $\tau$-admissibility and Property~2, the curve $\gamma$ is constant on $[a_{j-1}, a_j]$ and has zero $\tau$-length there). By Lemma~\ref{lem:2} and the left-invariance of $\tau$, we have $L_\tau(\gamma) = \sum_{j=1}^N L_\tau(\gamma|_{[a_{j-1}, a_j]}) = \sum_{j=1}^N L_\tau(\gamma_j)$, where $\gamma_j(s) := \gamma(a_{j-1})^{-1}\gamma(s)$ ($s \in [a_{j-1}, a_j]$). If we show that each non-constant arc $\gamma_j$, after a suitable monotonic reparameterization, becomes a Lipschitz $C^+$-admissible curve $\tilde{\gamma}_j$ for which $L_\tau(\gamma_j) \le \ell(\tilde{\gamma}_j)$, then the concatenation of the curves $\gamma(a_{j-1})\tilde{\gamma}_j$ (followed by a linear rescaling of the domain to the interval $[0, 1]$) will yield a Lipschitz $C^+$-admissible curve $\tilde{\gamma} \in \Omega_{C^+}(g, h)$ with length $\ell(\tilde{\gamma}) = \sum_{j=1}^N \ell(\tilde{\gamma}_j) \ge \sum_{j=1}^N L_\tau(\gamma_j) = L_\tau(\gamma)$. By the definition of $\tau_c(g, h)$, this implies $L_\tau(\gamma) \le \ell(\tilde{\gamma}) \le \tau_c(g, h)$, and given the arbitrariness of $\gamma \in \Omega_\tau(g, h)$, it yields the required inequality $\tau(g, h) \le \tau_c(g, h)$.

Therefore, without loss of generality, we may assume that $g = e \neq h$, and the entire curve $\gamma : [0, 1] \to B_{\rho_0}^G(e)$ connects $\gamma(0) = e$ to $\gamma(1) = h \neq e$ and lies entirely in the ball $B_{\rho_0}^G(e)$. Denote $z(s) := \exp^{-1}(\gamma(s)) \in B_{\rho_0}(0)$ ($s \in [0, 1]$, $z(0) = 0$, $z(1) \neq 0$). For any $s_1, s_2 \in [0, 1]$ with $s_1 < s_2$, set $\tilde{x}(s_1, s_2) := \exp^{-1}(\gamma(s_1)^{-1}\gamma(s_2)) = (-z(s_1)) * z(s_2)$. By Lemma~\ref{lember:2}, we have $|\tilde{x}(s_1, s_2)| \le 2\rho_0 + \rho_0^2 < 3\rho_0 < \frac{r}{2^{k_2}}$, and moreover $\tau(e, \exp(\tilde{x}(s_1, s_2))) = \tau(\gamma(s_1), \gamma(s_2)) \ge 0$. Applying Lemma~\ref{lem:6} to the vector $\tilde{x} := \tilde{x}(s_1, s_2)$ and Lemma~\ref{lember:1} to the equality $z(s_2) = z(s_1) * \tilde{x} = z(s_1) + \tilde{x} + z(z(s_1), \tilde{x})$, taking into account $\rho_0 \le \frac{c_0 r}{4} < \frac{c_0}{2}$, we obtain:
\begin{equation}\label{eq:xi_monotone}
\langle \xi, z(s_2) - z(s_1) \rangle = \langle \xi, \tilde{x} \rangle + \langle \xi, z(z(s_1), \tilde{x}) \rangle \ge c_0 |\tilde{x}| - \rho_0 |\tilde{x}| \ge \frac{c_0}{2}|\tilde{x}|.
\end{equation}
On the other hand, $|z(s_2) - z(s_1)| \le |\tilde{x}| + |z(z(s_1), \tilde{x})| \le (1 + \rho_0)|\tilde{x}| < 2|\tilde{x}|$. Combining this with~\eqref{eq:xi_monotone}, we find:
\begin{equation}\label{eq:lip_reparam}
|z(s_2) - z(s_1)| \le 2|\tilde{x}(s_1, s_2)| \le \frac{4}{c_0}\langle \xi, z(s_2) - z(s_1) \rangle \qquad (0 \le s_1 < s_2 \le 1).
\end{equation}

Introduce the continuous non-decreasing function $\sigma(s) := \langle \xi, z(s) \rangle$ ($s \in [0, 1]$). Since $z(0) = 0$ and $z(1) \neq 0$, \eqref{eq:lip_reparam} implies that $\sigma(0) = 0$ and $T := \sigma(1) \ge \frac{c_0}{4}|z(1)| > 0$; by the continuity of $\sigma$, its image is the entire interval $[0, T]$. Furthermore, by the estimate~\eqref{eq:lip_reparam}, the equality $\sigma(s_1) = \sigma(s_2)$ implies $z(s_1) = z(s_2)$. Consequently, the formula $\tilde{z}(t) := z(s)$ for $t = \sigma(s) \in [0, T]$ correctly defines a mapping $\tilde{z} : [0, T] \to B_{\rho_0}(0)$ (so that $z(s) = \tilde{z}(\sigma(s))$ for all $s \in [0, 1]$), which, due to~\eqref{eq:lip_reparam}, satisfies the Lipschitz condition with the constant $\frac{4}{c_0}$:
\[
|\tilde{z}(t_2) - \tilde{z}(t_1)| \le 2\big|\exp^{-1}\!\big(\tilde{\gamma}(t_1)^{-1}\tilde{\gamma}(t_2)\big)\big| \le \frac{4}{c_0}(t_2 - t_1) \qquad (0 \le t_1 < t_2 \le T),
\]
where $\tilde{\gamma}(t) := \exp(\tilde{z}(t))$ is a Lipschitz reparameterization of the curve $\gamma$ that preserves its $\tau$-length: $L_\tau(\tilde{\gamma}) = L_\tau(\gamma)$.

Consider the non-decreasing function $F(t) := L_\tau(\tilde{\gamma}|_{[0, t]})$ on the interval $[0, T]$. By Lemma~\ref{lem:2} and the second inequality~\eqref{eq:local_cone_bound} from Lemma~\ref{lem:6}, for any $0 \le t_1 < t_2 \le T$, we have:
\[
0 \le F(t_2) - F(t_1) = L_\tau(\tilde{\gamma}|_{[t_1, t_2]}) \le \tau(\tilde{\gamma}(t_1), \tilde{\gamma}(t_2)) \le \frac{2\tau(e, q)}{r}\big|\exp^{-1}\!\big(\tilde{\gamma}(t_1)^{-1}\tilde{\gamma}(t_2)\big)\big| \le \frac{4\tau(e, q)}{c_0 r}(t_2 - t_1).
\]
This means the function $F$ is Lipschitz on $[0, T]$ and, consequently, absolutely continuous: $L_\tau(\gamma) = F(T) = \int_0^T F'(t)\,dt$. By Rademacher's theorem, the mapping $\tilde{z}$ (and hence the curve $\tilde{\gamma}$) and the function $F$ are differentiable for almost all $t \in [0, T]$. Fix a point $t \in (0, T)$ at which the derivatives $\dot{\tilde{\gamma}}(t)$ and $F'(t)$ exist, and denote $u(t) := \theta\big(\dot{\tilde{\gamma}}(t)\big) = (dl_{\tilde{\gamma}(t)})_e^{-1}\big(\dot{\tilde{\gamma}}(t)\big) \in L$. By the definition of the derivative, we have:
\[
v_k := 2^k \exp^{-1}\!\big(\tilde{\gamma}(t)^{-1}\tilde{\gamma}(t + 2^{-k})\big) \to u(t) \qquad \text{as } k \to \infty.
\]

Since the curve $\tilde{\gamma}$ is $\tau$-admissible, for all sufficiently large $k$, we have $\tau(e, \exp(2^{-k}v_k)) = \tau(\tilde{\gamma}(t), \tilde{\gamma}(t + 2^{-k})) \ge 0$. Fix an arbitrary sufficiently large integer $m \ge 1$ (such that $\frac{|u(t)|}{2^m} < \lambda$, and therefore $\frac{v_k}{2^m} \in V$ for all large $k$); then for $k \ge m$, by the reverse triangle inequality, we have
\[
2^m \tau\!\left(e, \exp\!\left(\frac{v_k}{2^m}\right)\right) \ge 2^k \tau\!\left(e, \exp\!\left(\frac{v_k}{2^k}\right)\right) = 2^k \tau(\tilde{\gamma}(t), \tilde{\gamma}(t + 2^{-k})) \ge \frac{F(t + 2^{-k}) - F(t)}{2^{-k}} \ge 0.
\]
Taking the limit superior as $k \to \infty$ in this chain of inequalities and using the upper semicontinuity of $\tau$ on $U \times U \subset U_0 \times U_0$, we obtain
\[
2^m \tau\!\left(e, \exp\!\left(\frac{u(t)}{2^m}\right)\right) \ge F'(t) \ge 0
\]
for all sufficiently large $m \ge 1$. By Lemma~\ref{lem:5}, this implies that $u(t) \in C_0^+$, and passing to the limit as $m \to \infty$ yields $\nu_0(u(t)) \ge F'(t)$. Thus, the curve $\tilde{\gamma}$ (after a linear change of parameter $t = T s$, $s \in [0, 1]$) belongs to $\Omega_{C^+}(e, h)$, and moreover
\[
L_\tau(\gamma) = F(T) = \int_0^T F'(t)\,dt \le \int_0^T \nu_0(u(t))\,dt = \ell(\tilde{\gamma}).
\]
Consequently, $\tau(g, h) \le \tau_c(g, h)$ for all $g, h \in G$.

Let us now prove the reverse inequality $\tau_c(g, h) \le \tau(g, h)$ for all $g, h \in G$. If $\Omega_{C^+}(g, h) = \varnothing$, then by definition $\tau_c(g, h) = -\infty \le \tau(g, h)$. Suppose $\Omega_{C^+}(g, h) \neq \varnothing$. Note that for any absolutely continuous curve $\gamma \in \Omega_{C^+}(g, h)$ with an integrable left-invariant velocity $v(s) := (dl_{\gamma(s)})_e^{-1}(\dot{\gamma}(s)) \in L^1([0, 1]; L)$, the monotonic absolutely continuous change of parameter
\[
t(s) := \frac{1}{S}\int_0^s \big(1 + |v(\zeta)|\big)\,d\zeta, \qquad \text{where } S := \int_0^1 \big(1 + |v(\zeta)|\big)\,d\zeta \ge 1,
\]
transforms $\gamma$ into a Lipschitz curve $X(t) := \gamma(s(t)) \in \Omega_{C^+}(g, h)$ with a bounded left-invariant velocity $|u(t)| = \frac{S|v(s(t))|}{1 + |v(s(t))|} \le S$ for a.e. $t \in [0, 1]$ and the same sub-Lorentzian length $\ell(X) = \ell(\gamma)$. Therefore, it suffices to consider an arbitrary Lipschitz $C^+$-admissible curve $X \in \Omega_{C^+}(g, h)$ connecting $X(0) = g$ to $X(1) = h$ and to show that the curve $X$ is $\tau$-admissible and satisfies the inequality $\ell(X) \le L_\tau(X) \le \tau(g, h)$.

Due to the continuity of the curve $X$ on the interval $[0, 1]$, we can choose a partition $0 = a_0 < a_1 < \dots < a_N = 1$ fine enough so that for each $j = 1, \dots, N$ and all $t \in [a_{j-1}, a_j]$, the inclusion $X(a_{j-1})^{-1}X(t) \in B_{\rho_0}^G(e) \subset B_{\frac{r}{2}}^G(e)$ holds. If we prove that each arc $X_j(t) := X(a_{j-1})^{-1}X(t)$ ($t \in [a_{j-1}, a_j]$) is $\tau$-admissible and satisfies the estimate $\ell(X_j) \le L_\tau(X_j)$, then by the left-invariance and the reverse triangle inequality for $\tau$, the entire curve $X$ will also be $\tau$-admissible (meaning $X \in \Omega_\tau(g, h)$ and $\tau(g, h) \ge 0$). Furthermore, by the additivity of the functionals $\ell$ and $L_\tau$ (Lemma~\ref{lem:2}), we will have
\[
\ell(X) = \sum_{j=1}^N \ell(X_j) \le \sum_{j=1}^N L_\tau(X_j) = L_\tau(X) \le \tau(g, h),
\]
from which taking the supremum over all $X \in \Omega_{C^+}(g, h)$ yields the required inequality $\tau_c(g, h) \le \tau(g, h)$.

Thus, applying a linear change of parameter on the interval $[0, 1]$ if necessary, we may assume without loss of generality that $g = e$ and the entire Lipschitz $C^+$-admissible curve $X : [0, 1] \to B_{\rho_0}^G(e)$ starts at $X(0) = e$ and lies entirely in the ball $B_{\rho_0}^G(e) \subset B_{\frac{r}{2}}^G(e)$. At the same time, for any subinterval $[t_1, t_2] \subset [0, 1]$, the curve $\widetilde{X}(\theta) := X(t_1)^{-1}X\big(t_1 + \theta(t_2 - t_1)\big)$, $\theta \in [0, 1]$, obtained by left translation by $X(t_1)^{-1}$ and a linear change of parameter, starts at $\widetilde{X}(0) = e$ and, by Lemma~\ref{lember:2}, lies entirely in the ball $B_{3\rho_0}^G(e) \subset B_{\frac{r}{2}}^G(e)$ (since $3\rho_0 < \frac{3}{8}r < \frac{r}{2}$). Therefore, to establish the inequality $\tau(X(t_1), X(t_2)) \ge \int_{t_1}^{t_2} \nu_0(u(s))\,ds \ge 0$ for all $0 \le t_1 < t_2 \le 1$, it suffices to prove the estimate $\tau(e, X(1)) \ge \ell(X)$ for an arbitrary Lipschitz $C^+$-admissible curve $X : [0, 1] \to B_{\frac{r}{2}}^G(e)$ with $X(0) = e$.

Denote by $u(t) := (dl_{X(t)})_e^{-1}\big(\dot{X}(t)\big) \in C_0^+$ ($t \in [0, 1]$) the measurable essentially bounded control of the curve $X$: there exists a constant $K > 0$ such that $|u(t)| \le K$ for almost all $t \in [0, 1]$. Set $z(t) := \exp^{-1}(X(t)) \in B_{\frac{r}{2}}(0)$ ($z(0) = 0$). In exponential coordinates, the curve $z(t)$ satisfies the differential equation
\begin{equation}\label{eq:ode_z}
\dot{z}(t) = \Psi(z(t))u(t) \qquad \text{for a.e. } t \in [0, 1],
\end{equation}
where $\Psi(z) := (d\exp_z)^{-1} \circ (dl_{\exp(z)})_e : L \to L$ is a smooth family of linear operators for $z \in B_r(0)$, with $\Psi(0) = \mathrm{id}_L$. By the compactness of the ball $B_r(0)$, there exist constants $M_1, M_2 > 0$ such that for all $z, w \in B_r(0)$, the following operator norm estimates hold:
\begin{equation}\label{eq:Psi_bounds}
\|\Psi(z)\| \le M_1, \qquad \|\Psi(z) - \Psi(w)\| \le M_2 |z - w|.
\end{equation}
For each integer $m \ge 1$, we partition the interval $[0, 1]$ by the points $t_i := \frac{i}{m}$ ($i = 0, 1, \dots, m$) into $m$ equal segments of length $\Delta t := \frac{1}{m}$, and define the average values of the control on each segment:
\[
u_{m, i} := \frac{1}{\Delta t}\int_{t_{i-1}}^{t_i} u(s)\,ds \qquad (i = 1, \dots, m).
\]
Since the cone $C_0^+$ is convex and closed, and $|u(s)| \le K$ a.e., we have $u_{m, i} \in C_0^+$ and $|u_{m, i}| \le K$. We will show that for the function $\nu_0$, the integral Jensen's inequality holds on each segment $[t_{i-1}, t_i]$:
\begin{equation}\label{eq:jensen_nu0}
\Delta t\,\nu_0(u_{m, i}) \ge \int_{t_{i-1}}^{t_i} \nu_0(u(s))\,ds.
\end{equation}
Fix a vector $u_0 \in \operatorname{int}_{L_0}(C_0^+)$. For any $\delta > 0$ and any $w \in C_0^+$, the vector $w + \delta u_0$ belongs to the relative interior $\operatorname{int}_{L_0}(C_0^+)$ of the cone $C_0^+$ in the subspace $L_0 = C_0^+ - C_0^+$. In particular, the point $v_\delta := u_{m, i} + \delta u_0$ lies in $\operatorname{int}_{L_0}(C_0^+)$. Since the function $-\nu_0$ is convex and finite on $C_0^+$, by \cite[Theorem~23.4]{Rockafellar1970} it has a subgradient $-\phi_\delta \in L_0$ at the interior point $v_\delta \in \operatorname{int}_{L_0}(C_0^+)$ (that is, $\phi_\delta \in L_0$ is a supergradient of the concave function $\nu_0$ at the point $v_\delta$), such that
\[
\nu_0(y) \le \nu_0(v_\delta) + \langle \phi_\delta, y - v_\delta \rangle \qquad \forall y \in C_0^+.
\]
Substituting $y = u(s) + \delta u_0 \in C_0^+$ and noting that $\nu_0(u(s) + \delta u_0) \ge \nu_0(u(s)) + \delta \nu_0(u_0) \ge \nu_0(u(s))$, we obtain
\[
\nu_0(u(s)) \le \nu_0(u_{m, i} + \delta u_0) + \langle \phi_\delta, u(s) - u_{m, i} \rangle \qquad \text{for a.e. } s \in [t_{i-1}, t_i].
\]
Integrating this inequality over $s \in [t_{i-1}, t_i]$ and observing that $\int_{t_{i-1}}^{t_i} \langle \phi_\delta, u(s) - u_{m, i} \rangle\,ds = 0$ by the definition of the average $u_{m, i}$, we find
\[
\int_{t_{i-1}}^{t_i} \nu_0(u(s))\,ds \le \Delta t\,\nu_0(u_{m, i} + \delta u_0).
\]
Taking the limit superior as $\delta \to 0^+$ and using the upper semicontinuity of the function $\nu_0$ on $C_0^+$, we obtain inequality~\eqref{eq:jensen_nu0}.

Now we construct a piecewise-exponential curve $X_m : [0, 1] \to G$ by setting $X_m(0) := e$ and
\[
X_m(t) := X_m(t_{i-1})\exp\big((t - t_{i-1})u_{m, i}\big) \qquad \text{for } t \in [t_{i-1}, t_i],\quad i = 1, \dots, m.
\]
We will show that for all sufficiently large $m$, the curve $X_m([0, 1])$ is entirely contained in the ball $B_r^G(e)$, and its endpoint $X_m(1)$ converges to $X(1)$ as $m \to \infty$. Let $[0, t_m^*] \subset [0, 1]$ be the maximal interval on which $X_m(t) \in B_r^G(e)$, and let $z_m(t) := \exp^{-1}(X_m(t)) \in B_r(0)$ for $t \in [0, t_m^*]$. On each segment $t \in [t_{i-1}, t_i] \cap [0, t_m^*]$, the curve $X_m(t)$ has a left-invariant velocity $u_m(t) := u_{m, i}$; hence $z_m(t)$ satisfies the equation $\dot{z}_m(t) = \Psi(z_m(t))u_m(t)$. Subtracting this from \eqref{eq:ode_z} and integrating from $0$ to $t \in [t_{k-1}, t_k] \cap [0, t_m^*]$, we obtain:
\begin{align*}
z(t) - z_m(t) &= \int_0^t \big(\Psi(z(s)) - \Psi(z_m(s))\big)u(s)\,ds + \int_0^t \Psi(z_m(s))\big(u(s) - u_m(s)\big)\,ds \\
&= \int_0^t \big(\Psi(z(s)) - \Psi(z_m(s))\big)u(s)\,ds + \sum_{i=1}^{k-1} \int_{t_{i-1}}^{t_i} \big(\Psi(z_m(s)) - \Psi(z_m(t_{i-1}))\big)\big(u(s) - u_{m, i}\big)\,ds \\
&\quad + \int_{t_{k-1}}^t \Psi(z_m(s))\big(u(s) - u_{m, k}\big)\,ds,
\end{align*}
since $\int_{t_{i-1}}^{t_i} \Psi(z_m(t_{i-1}))(u(s) - u_{m, i})\,ds = 0$. In view of the estimates~\eqref{eq:Psi_bounds} and $|\dot{z}_m(s)| \le M_1 K$ for $s \in [t_{i-1}, t_i] \cap [0, t_m^*]$, we have $|z_m(s) - z_m(t_{i-1})| \le M_1 K \Delta t = \frac{M_1 K}{m}$, which yields $\|\Psi(z_m(s)) - \Psi(z_m(t_{i-1}))\| \le \frac{M_1 M_2 K}{m}$. Consequently,
\[
|z(t) - z_m(t)| \le M_2 K \int_0^t |z(s) - z_m(s)|\,ds + \frac{2M_1 M_2 K^2}{m} + \frac{2M_1 K}{m}.
\]
By Gronwall's lemma, for all $t \in [0, t_m^*]$ the following estimate holds:
\[
|z(t) - z_m(t)| \le \frac{2M_1 K(M_2 K + 1)}{m} e^{M_2 K} \to 0 \qquad \text{as } m \to \infty.
\]
Since $|z(t)| \le \frac{r}{2}$, for all sufficiently large $m$ the right-hand side is less than $\frac{r}{4}$, so $|z_m(t)| \le |z(t)| + |z(t) - z_m(t)| < \frac{3}{4}r < r$ on $[0, t_m^*]$. By the maximality of the interval $[0, t_m^*]$, it follows that $t_m^* = 1$, meaning $|z_m(t)| < r$ for all $t \in [0, 1]$, and $X_m(1) = \exp(z_m(1)) \to \exp(z(1)) = X(1)$ as $m \to \infty$.

Finally, we establish a lower bound for $\tau(e, X_m(1))$. By Lemma~\ref{lem:5}, for each vector $u_{m, i} \in C_0^+$ we have $\tau(e, \exp(\Delta t\,u_{m, i})) \ge \Delta t\,\nu_0(u_{m, i})$. Applying the reverse triangle inequality, the left-invariance of $\tau$, and Jensen's inequality \eqref{eq:jensen_nu0}, we find:
\begin{align*}
\tau(e, X_m(1)) 
&\ge \sum_{i=1}^m \tau\big(X_m(t_{i-1}), X_m(t_i)\big) \\
&= \sum_{i=1}^m \tau\big(e, \exp(\Delta t\,u_{m, i})\big) \\
&\ge \sum_{i=1}^m \Delta t\,\nu_0(u_{m, i}) \\
&\ge \int_0^1 \nu_0(u(s))\,ds = \ell(X).
\end{align*}
Taking the limit superior as $m \to \infty$ in the inequality $\tau(e, X_m(1)) \ge \ell(X)$ and using the upper semicontinuity of $\tau$ on $U \times U \subset U_0 \times U_0$, we obtain $\tau(e, X(1)) \ge \ell(X) \ge 0$. Applying the exact same reasoning to the curve $\widetilde{X}$ for any subinterval $[t_1, t_2] \subset [0, 1]$, we conclude that $\tau(X(t_1), X(t_2)) \ge \ell(X|_{[t_1, t_2]}) \ge 0$, meaning $X$ is a $\tau$-admissible curve. Thus, for any partition $0 = s_0 < s_1 < \dots < s_p = 1$, we have
\[
\sum_{k=1}^p \tau\big(X(s_{k-1}), X(s_k)\big) \ge \sum_{k=1}^p \ell\big(X|_{[s_{k-1}, s_k]}\big) = \ell(X).
\]
Taking the infimum over all partitions of the interval $[0, 1]$, we obtain $L_\tau(X) \ge \ell(X)$, whence $\tau(g, h) \ge \tau_c(g, h)$. Thus, $\tau(g, h) = \tau_c(g, h)$ for all $g, h \in G$.

Let us now prove the uniqueness of the pair $(C_0^+, \nu_0)$. Suppose $(\widetilde{C}_0^+, \widetilde{\nu}_0)$ is another pair defining a left-invariant sub-Lorentzian structure on $G$, whose generated chronometer $\widetilde{\tau}_c$ coincides with $\tau$. For any $u \in \widetilde{C}_0^+$ and any $t > 0$, the curve $\gamma_t(s) := \exp(stu)$ ($s \in [0, 1]$) belongs to $\Omega_{\widetilde{C}^+}(e, \exp(tu))$, therefore
\[
\tau(e, \exp(tu)) = \widetilde{\tau}_c(e, \exp(tu)) \ge \int_0^1 \widetilde{\nu}_0(tu)\,ds = t\,\widetilde{\nu}_0(u) \ge 0.
\]
By Lemma~\ref{lem:5}, this implies that $u \in C_0^+$ and $\nu_0(u) = \lim_{t \to 0^+} \frac{\tau(e, \exp(tu))}{t} \ge \widetilde{\nu}_0(u)$. Thus, $\widetilde{C}_0^+ \subset C_0^+$ and $\widetilde{\nu}_0 \le \nu_0$ on $\widetilde{C}_0^+$.

Conversely, let $u \in C_0^+ \setminus \{0\}$. Applying Lemma~\ref{lem:1} to the neighborhood $\operatorname{int} B_{\rho_0}^G(e)$, we see that for all sufficiently large integers $m$, the point $g_m := \exp\!\left(\frac{u}{m}\right) \to e$ satisfies $J(e, g_m) \subset B_{\rho_0}^G(e)$, and by Lemma~\ref{lem:5} we have $\widetilde{\tau}_c(e, g_m) = \tau(e, g_m) \ge \frac{1}{m}\nu_0(u) \ge 0$. We can choose a curve $X_m \in \Omega_{\widetilde{C}^+}(e, g_m)$ with control $\tilde{u}_m(s) := (dl_{X_m(s)})_e^{-1}\big(\dot{X}_m(s)\big) \in \widetilde{C}_0^+ \subset C_0^+$ such that $\int_0^1 \widetilde{\nu}_0(\tilde{u}_m(s))\,ds \ge \frac{1}{m}\nu_0(u) - \frac{1}{m^2}$. Since $\tilde{u}_m(s) \in C_0^+$, the curve $X_m$ is $C^+$-admissible and, as proven above, $\tau$-admissible, whence $X_m([0, 1]) \subset J(e, g_m) \subset B_{\rho_0}^G(e)$. Setting $z_m(s) := \exp^{-1}(X_m(s))$ ($z_m(0) = 0$, $z_m(1) = \frac{u}{m}$), from the estimates~\eqref{eq:lip_reparam} and~\eqref{eq:Psi_bounds} we obtain $|z_m(s)| \le \frac{4}{c_0}\langle \xi, z_m(1) \rangle \le \frac{4|u|}{c_0 m}$ and $\|\Psi(z_m(s)) - \mathrm{id}_L\| \le \frac{4M_2 |u|}{c_0 m} \le c_0$ (for large $m$). Then, by part 1 of Lemma~\ref{lem:3}, the equation $\dot{z}_m(s) = \Psi(z_m(s))\tilde{u}_m(s)$ implies $\langle \xi, \dot{z}_m(s) \rangle \ge 2c_0 |\tilde{u}_m(s)| - c_0 |\tilde{u}_m(s)| = c_0 |\tilde{u}_m(s)|$, and integrating over $s \in [0, 1]$ yields $\int_0^1 |\tilde{u}_m(s)|\,ds \le \frac{1}{c_0}\langle \xi, z_m(1) \rangle \le \frac{|u|}{c_0 m}$. Therefore, for the vector $v_m := m \int_0^1 \tilde{u}_m(s)\,ds \in \widetilde{C}_0^+$ we have
\[
|u - v_m| = m\left|\int_0^1 \big(\Psi(z_m(s)) - \mathrm{id}_L\big)\tilde{u}_m(s)\,ds\right| \le m \cdot \frac{4M_2 |u|}{c_0 m} \cdot \frac{|u|}{c_0 m} = \frac{4M_2 |u|^2}{c_0^2 m} \to 0 \qquad (m \to \infty).
\]
Due to the closedness of the cone $\widetilde{C}_0^+$, this implies $u = \lim_{m \to \infty} v_m \in \widetilde{C}_0^+$. By positive homogeneity, Jensen's inequality~\eqref{eq:jensen_nu0} (applied to the anti-norm $\widetilde{\nu}_0$), and the upper semicontinuity of $\widetilde{\nu}_0$, we find
\[
\widetilde{\nu}_0(u) \ge \varlimsup_{m \to \infty} \widetilde{\nu}_0(v_m) \ge \varlimsup_{m \to \infty} m \int_0^1 \widetilde{\nu}_0(\tilde{u}_m(s))\,ds \ge \nu_0(u).
\]
Consequently, $\widetilde{C}_0^+ = C_0^+$ and $\widetilde{\nu}_0 = \nu_0$. This completes the proof of the theorem.
\end{proof}

\section{Direct Theorem}

\begin{definition}[Chronological condition]\label{def:condition}
We say that a cone $C_0^+ \subset L$ satisfies the chronological condition if for any neighborhood $U$ of the identity $e \in G$, there exists a neighborhood $V \subset U$ of the identity $e$ such that every $C^+$-admissible curve $\gamma : [0, 1] \to G$ with endpoints $\gamma(0), \gamma(1) \in V$ is entirely contained in $U$:
\[
\gamma(0), \gamma(1) \in V \implies \gamma([0, 1]) \subset U
\]
(in other words, $J_{C^+}(V, V) \subset U$), where
\begin{align*}
J_{C^+}(A, B) := \big\{ x \in G : \exists\, p \in A,\, q \in B \text{ such that } & \Omega_{C^+}(p, x) \neq \varnothing, \\
& \Omega_{C^+}(x, q) \neq \varnothing \big\}.
\end{align*}
\end{definition}

\begin{remark}
Note that the cone $C_0^+$ constructed in Theorem~\ref{thm:main} from the chronometer $\tau$ automatically satisfies the chronological condition. Indeed, for any neighborhood $U$ of the identity $e \in G$, by Property~5 there exist points $p_0, q_0 \in U$ such that $e \in V := \operatorname{int} I(p_0, q_0) \subset I(p_0, q_0) \subset U$. If $x \in J_{C^+}(V, V)$, then there exist points $p, q \in V \subset I(p_0, q_0)$ for which $\Omega_{C^+}(p, x) \neq \varnothing$ and $\Omega_{C^+}(x, q) \neq \varnothing$. Then by Theorem~\ref{thm:main}, we have $\tau(p, x) = \tau_c(p, x) \ge 0$ and $\tau(x, q) = \tau_c(x, q) \ge 0$, whence by the reverse triangle inequality (Property~3) we obtain $\tau(p_0, x) \ge \tau(p_0, p) + \tau(p, x) > 0$ and $\tau(x, q_0) \ge \tau(x, q) + \tau(q, q_0) > 0$. Consequently, $x \in I(p_0, q_0) \subset U$, meaning $J_{C^+}(V, V) \subset U$.
\end{remark}

\begin{theorem}\label{thm:direct}
Let $(G, C_0^+, \nu_0)$ be a left-invariant sub-Lorentzian structure on a connected finite-dimensional Lie group $G$ with Lie algebra $L$; that is, $C_0^+ \subset L$ is a non-trivial closed pointed convex cone whose linear span $L_0 = C_0^+ - C_0^+$ generates the entire Lie algebra $L$ ($\widetilde{L}_0 = L$), and $\nu_0 : L \to [0, +\infty) \cup \{-\infty\}$ is an anti-norm on the cone $C_0^+$ ($\nu_0|_{C_0^+} \not\equiv 0$). Suppose that the cone $C_0^+$ satisfies the chronological condition.
Then the sub-Lorentzian chronometer generated by this structure,
\[
\tau_c(g, h) := \sup_{\gamma \in \Omega_{C^+}(g, h)} \int_0^1 \nu_0\big(\theta(\dot{\gamma}(s))\big)\,ds \qquad (\sup \varnothing := -\infty),
\]
satisfies Properties~1--6.
\end{theorem}

\begin{proof}[Proof of Theorem \ref{thm:direct}]
\textit{Preliminary local estimates.} Since $C_0^+ \subset L$ is a non-trivial closed pointed convex cone, just as in part 1 of Lemma~\ref{lem:3}, the Hahn-Banach separation theorem implies the existence of a unit vector $\xi \in L$ ($|\xi| = 1$) and a constant $c_0 \in (0, 1/2)$ such that
\begin{equation}\label{eq:direct_xi}
\langle \xi, u \rangle \ge 2c_0 |u| \qquad \forall u \in C_0^+.
\end{equation}
Furthermore, since the function $\nu_0$ is upper semicontinuous and finite on the compact set $C_0^+ \cap \partial B_1(0)$, it attains its supremum $M_\nu := \max\limits_{u \in C_0^+,\, |u| = 1} \nu_0(u) \in (0, +\infty)$ on it. Due to the positive homogeneity and non-negativity of $\nu_0$ on $C_0^+$, this implies the bound
\begin{equation}\label{eq:direct_nu_bound}
0 \le \nu_0(u) \le M_\nu |u| \qquad \forall u \in C_0^+.
\end{equation}
Let $B_r(0) \subset V$ be a closed ball in $L$ on which the exponential map $\exp$ is a diffeomorphism and the operator $\Psi(z) = (d\exp_z)^{-1} \circ (dl_{\exp(z)})_e$ from equation \eqref{eq:ode_z} satisfies the estimates \eqref{eq:Psi_bounds} with constants $M_1, M_2 > 0$. We choose a radius $r_0 \in (0, r]$ small enough so that in the exponential coordinates $z = \exp^{-1}(x) \in B_{r_0}(0)$, the inequality $\|\Psi(z) - \mathrm{id}_L\| \le M_2 r_0 \le c_0$ holds for all $z \in B_{r_0}(0)$. Then, for any $C^+$-admissible curve $\gamma : [0, 1] \to B_{r_0}^G(e)$ with control $u(s) := \theta(\dot{\gamma}(s)) \in C_0^+$, its coordinate representation $z(s) := \exp^{-1}(\gamma(s)) \in B_{r_0}(0)$ satisfies the equation $\dot{z}(s) = \Psi(z(s))u(s)$, whence by \eqref{eq:direct_xi}, for almost all $s \in [0, 1]$, we have:
\[
\langle \xi, \dot{z}(s) \rangle = \langle \xi, u(s) \rangle + \big\langle \xi, \big(\Psi(z(s)) - \mathrm{id}_L\big)u(s) \big\rangle \ge 2c_0 |u(s)| - c_0 |u(s)| = c_0 |u(s)|.
\]
Integrating this inequality over the interval $[0, 1]$ and taking \eqref{eq:direct_nu_bound} into account, we obtain the a priori estimates:
\begin{equation}\label{eq:direct_L1_est}
\int_0^1 |u(s)|\,ds \le \frac{1}{c_0}\langle \xi, z(1) - z(0) \rangle \le \frac{1}{c_0}|z(1) - z(0)|, \qquad \ell(\gamma) \le \frac{M_\nu}{c_0}|z(1) - z(0)|.
\end{equation}
By the chronological condition, for the neighborhood $B_{r_0}^G(e)$ of the identity $e$, there exists an open neighborhood $V_* \subset B_{r_0}^G(e)$ of $e$ such that any $C^+$-admissible curve with endpoints in $V_*$ is entirely contained in $B_{r_0}^G(e)$. We choose an open neighborhood $U_0 \subset V_*$ of $e$ such that $U_0^{-1}U_0 \subset V_*$.

Let us now sequentially verify that the function $\tau_c$ satisfies Properties 1--6.

\textit{Properties 1 and 3 (Left-invariance and reverse triangle inequality).} For any $a, g, h \in G$, the left translation $\gamma \mapsto l_a \circ \gamma$ establishes a bijection between the sets $\Omega_{C^+}(g, h)$ and $\Omega_{C^+}(ag, ah)$ that preserves the left-invariant velocity $\theta(\dot{\gamma}(s))$ and, consequently, the length functional $\ell(\gamma)$. Taking the supremum, we obtain $\tau_c(ag, ah) = \tau_c(g, h)$ (Property 1).

Suppose $\tau_c(g_1, g_2) \ge 0$ and $\tau_c(g_2, g_3) \ge 0$; then the sets $\Omega_{C^+}(g_1, g_2)$ and $\Omega_{C^+}(g_2, g_3)$ are non-empty. For any curves $\gamma_1 \in \Omega_{C^+}(g_1, g_2)$ and $\gamma_2 \in \Omega_{C^+}(g_2, g_3)$, their concatenation $\gamma(s)$, which equals $\gamma_1(2s)$ for $s \in [0, 1/2]$ and $\gamma_2(2s - 1)$ for $s \in [1/2, 1]$, belongs to $\Omega_{C^+}(g_1, g_3)$ and satisfies the equality $\ell(\gamma) = \ell(\gamma_1) + \ell(\gamma_2) \le \tau_c(g_1, g_3)$. Taking the supremum over all $\gamma_1$ and $\gamma_2$, we obtain $\tau_c(g_1, g_3) \ge \tau_c(g_1, g_2) + \tau_c(g_2, g_3)$ (Property 3).

\textit{Property 2 (Reflexivity and antisymmetry).} The constant curve $\gamma_0(s) \equiv e$ belongs to $\Omega_{C^+}(e, e)$ and has length $\ell(\gamma_0) = 0$, thus $\tau_c(e, e) \ge 0$. Conversely, for any curve $\gamma \in \Omega_{C^+}(e, e)$, we have $\gamma(0) = \gamma(1) = e \in V_*$, whence by the chronological condition, the entire image $\gamma([0, 1])$ lies in $B_{r_0}^G(e)$. Then the estimate \eqref{eq:direct_L1_est} for $z(0) = z(1) = 0$ implies $\int_0^1 |u(s)|\,ds = 0$, meaning $u(s) = 0$ a.e., $\gamma(s) \equiv e$, and $\ell(\gamma) = 0$. Hence, $\tau_c(g, g) = \tau_c(e, e) = 0$ for all $g \in G$.

If, on the other hand, $g \neq h$, $\tau_c(g, h) \ge 0$, and (by contradiction) $\tau_c(h, g) \ge 0$, then there exist curves $\gamma_1 \in \Omega_{C^+}(e, g^{-1}h)$ and $\gamma_2 \in \Omega_{C^+}(g^{-1}h, e)$. Their concatenation $\gamma$ belongs to $\Omega_{C^+}(e, e)$, yet passes through the point $\gamma(1/2) = g^{-1}h \neq e$, which contradicts the previously proven identity $\gamma(s) \equiv e$ for all $\gamma \in \Omega_{C^+}(e, e)$. Consequently, $\tau_c(h, g) = -\infty$.

\textit{Property 4 (Local upper semicontinuity).} We will show that for the previously chosen neighborhood $U_0$ and any $c \ge 0$, the set $E_c := \{(g, h) \in U_0 \times U_0 : \tau_c(g, h) \ge c\}$ is closed in $U_0 \times U_0$. Let a sequence $(g_n, h_n) \in E_c$ converge to $(g, h) \in U_0 \times U_0$. Set $x_n := g_n^{-1}h_n \in V_*$ and $x := g^{-1}h \in V_*$; then $x_n \to x$ and $\tau_c(e, x_n) = \tau_c(g_n, h_n) \ge c$.

For each $n \ge 1$, we choose a curve $\gamma_n \in \Omega_{C^+}(e, x_n)$ such that $\ell(\gamma_n) \ge c - \frac{1}{n}$ (note that $\tau_c(e, x_n) \le \frac{M_\nu r_0}{c_0} < +\infty$ by \eqref{eq:direct_L1_est}). Since $\gamma_n(0) = e \in V_*$ and $\gamma_n(1) = x_n \in V_*$, the choice of $V_*$ implies $\gamma_n([0, 1]) \subset B_{r_0}^G(e)$, and from the estimate \eqref{eq:direct_L1_est}, for the control $u_n(s) := \theta(\dot{\gamma}_n(s))$, it follows that $L_n := \int_0^1 |u_n(s)|\,ds \le \frac{r_0}{c_0}$. By performing a monotonic change of parameter $t(s) := \frac{1}{1 + L_n}\int_0^s (1 + |u_n(\zeta)|)\,d\zeta$ on $[0, 1]$, we transition to a curve $\tilde{\gamma}_n \in \Omega_{C^+}(e, x_n)$ with the same length $\ell(\tilde{\gamma}_n) = \ell(\gamma_n) \ge c - \frac{1}{n}$, whose control $\tilde{u}_n(t) \in C_0^+$ is uniformly bounded:
\[
|\tilde{u}_n(t)| \le 1 + L_n \le 1 + \frac{r_0}{c_0} =: K_0 \qquad \text{for a.e. } t \in [0, 1].
\]
Then the coordinate curves $\tilde{z}_n(t) := \exp^{-1}(\tilde{\gamma}_n(t)) \in B_{r_0}(0)$ are uniformly Lipschitz ($|\dot{\tilde{z}}_n(t)| \le M_1 K_0$ a.e.). By the Arzelà-Ascoli theorem and the weak compactness of closed bounded sets in the Hilbert space $L^2([0, 1]; L)$, passing to a subsequence, we may assume that $\tilde{z}_n \to \tilde{z}$ uniformly on $[0, 1]$, and $\tilde{u}_n \rightharpoonup \tilde{u}$ weakly in $L^2([0, 1]; L)$. Passing to the limit as $n \to \infty$ in the integral identity $\tilde{z}_n(t) = \int_0^t \Psi(\tilde{z}_n(s))\tilde{u}_n(s)\,ds$, we obtain $\tilde{z}(t) = \int_0^t \Psi(\tilde{z}(s))\tilde{u}(s)\,ds$, with $\tilde{z}(0) = 0$ and $\tilde{z}(1) = \lim_{n \to \infty} \tilde{z}_n(1) = \exp^{-1}(x)$.

By Mazur's lemma, for each $m \ge 1$ there exist convex combinations $w_m = \sum_{k=m}^{N_m} \alpha_{m, k}\tilde{u}_k$ ($\alpha_{m, k} \ge 0$, $\sum_{k=m}^{N_m} \alpha_{m, k} = 1$) that converge to $\tilde{u}$ strongly in $L^2([0, 1]; L)$, and after passing to a subsequence (for which we retain the previous notation), they converge pointwise a.e. on $[0, 1]$. Due to the convexity and closedness of $C_0^+$, we have $\tilde{u}(t) \in C_0^+$ a.e. Since $\dot{\tilde{z}}(t) = \Psi(\tilde{z}(t))\tilde{u}(t)$ a.e., the curve $\tilde{\gamma}(t) := \exp(\tilde{z}(t))$ has the control $\theta(\dot{\tilde{\gamma}}(t)) = \tilde{u}(t) \in C_0^+$ a.e., and therefore belongs to $\Omega_{C^+}(e, x)$. Since the function $\nu_0$ is concave and upper semicontinuous on $C_0^+$, and $\nu_0(w_m(t)) \le M_\nu K_0$, by Fatou's lemma we find:
\begin{align*}
\tau_c(g, h) = \tau_c(e, x) 
&\ge \int_0^1 \nu_0(\tilde{u}(t))\,dt \ge \int_0^1 \varlimsup_{m \to \infty} \nu_0(w_m(t))\,dt \ge \\
&\ge \varlimsup_{m \to \infty} \int_0^1 \nu_0(w_m(t))\,dt \ge \varlimsup_{m \to \infty} \sum_{k=m}^{N_m} \alpha_{m, k}\ell(\tilde{\gamma}_k) \ge c.
\end{align*}
Consequently, $(g, h) \in E_c$, which proves Property 4.

\textit{Property 5 (Compatibility with the manifold topology).} Due to left-invariance ($l_x(I(p, q)) = I(xp, xq)$), it suffices to verify the local condition of Property 5 at the point $x = e$ (from which it will immediately follow that the family $\{\operatorname{int} I(p, q)\}_{p, q \in G}$ forms a basis for the topology of the group $G$). Since the cone $C_0^+$ has a non-empty relative interior $\operatorname{int}_{L_0}(C_0^+)$ in its linear span $L_0 = C_0^+ - C_0^+$, and the Lie subalgebra $\widetilde{L}_0$ generated by $L_0$ coincides with the entire Lie algebra $L$, the Sussmann-Jurdjevic-Krener theorem on reachability guarantees that for any neighborhood $W$ of the identity $e$, the set of points of the form $y = \exp(t_1 v_1)\dots\exp(t_k v_k) \in W$, where $v_i \in \operatorname{int}_{L_0}(C_0^+)$ and $t_i > 0$, contains a non-empty open subset $W_0 \subset W$. Since $\nu_0(v_i) > 0$ for $v_i \in \operatorname{int}_{L_0}(C_0^+)$, for each point $y \in W_0$ the piecewise-exponential curve from $e$ to $y$ yields the estimate $\tau_c(e, y) \ge \sum_{i=1}^k t_i \nu_0(v_i) > 0$. Hence, $W_0 \subset I^+(e) \cap W$, and because $W_0$ is open, we have $W_0 \subset \operatorname{int} I^+(e) \cap W$, which implies $e \in \overline{\operatorname{int} I^+(e)}$. Because the mapping $y \mapsto y^{-1}$ is a diffeomorphism and $\tau_c(y^{-1}, e) = \tau_c(e, y)$, we similarly have $e \in \overline{\operatorname{int} I^-(e)}$.

Now let $U \subset G$ be an arbitrary neighborhood of the identity $e$. By the chronological condition, we can choose an open neighborhood $V \subset U$ of $e$ such that $J_{C^+}(V, V) \subset U$. Since $e \in \overline{\operatorname{int} I^-(e)}$ and $e \in \overline{\operatorname{int} I^+(e)}$, there exist points $p \in V \cap \operatorname{int} I^-(e)$ and $q \in V \cap \operatorname{int} I^+(e)$. Because the mappings $y \mapsto y^{-1}p$ and $y \mapsto y^{-1}q$ are homeomorphisms of the group $G$ onto itself mapping the identity $e$ to $p$ and to $q$, respectively, the left-invariance relations $\tau_c(p, y) = \tau_c(y^{-1}p, e)$ and $\tau_c(y, q) = \tau_c(e, y^{-1}q)$ imply that the conditions $p \in \operatorname{int} I^-(e)$ and $q \in \operatorname{int} I^+(e)$ are equivalent to the inclusions $e \in \operatorname{int} I^+(p)$ and $e \in \operatorname{int} I^-(q)$. Therefore,
\[
e \in \operatorname{int} I^+(p) \cap \operatorname{int} I^-(q) = \operatorname{int} I(p, q) \subset I(p, q).
\]
Moreover, for any point $y \in I(p, q)$, we have $\tau_c(p, y) > 0$ and $\tau_c(y, q) > 0$, meaning $\Omega_{C^+}(p, y) \neq \varnothing$ and $\Omega_{C^+}(y, q) \neq \varnothing$, whence $y \in J_{C^+}(V, V) \subset U$. Thus, $p, q \in V \subset U$ and $e \in \operatorname{int} I(p, q) \subset I(p, q) \subset U$.

\textit{Property 6 (Intrinsic chronometer).} Suppose $g \neq h$ and $\tau_c(g, h) \ge 0$; then $\Omega_{C^+}(g, h) \neq \varnothing$. For any curve $\gamma \in \Omega_{C^+}(g, h)$ and any $0 \le s_1 < s_2 \le 1$, the restriction $\gamma|_{[s_1, s_2]}$, after a linear change of parameter, belongs to $\Omega_{C^+}(\gamma(s_1), \gamma(s_2))$, so $\tau_c(\gamma(s_1), \gamma(s_2)) \ge \ell(\gamma|_{[s_1, s_2]}) \ge 0$. Consequently, $\gamma \in \Omega_{\tau_c}(g, h)$, meaning $\Omega_{C^+}(g, h) \subset \Omega_{\tau_c}(g, h)$, and the set $\Omega_{\tau_c}(g, h)$ is non-empty. Moreover, for any partition $0 = s_0 < s_1 < \dots < s_m = 1$, we have
\[
\sum_{i=1}^m \tau_c\big(\gamma(s_{i-1}), \gamma(s_i)\big) \ge \sum_{i=1}^m \ell\big(\gamma|_{[s_{i-1}, s_i]}\big) = \ell(\gamma),
\]
whence, by taking the infimum over all partitions, we obtain $L_{\tau_c}(\gamma) \ge \ell(\gamma)$ for all $\gamma \in \Omega_{C^+}(g, h)$. At the same time, for any continuous $\tau_c$-admissible curve $\tilde{\gamma} \in \Omega_{\tau_c}(g, h)$, by substituting the one-step partition $0 = s_0 < s_1 = 1$ into the definition of $L_{\tau_c}(\tilde{\gamma})$, we have $L_{\tau_c}(\tilde{\gamma}) \le \tau_c(\tilde{\gamma}(0), \tilde{\gamma}(1)) = \tau_c(g, h)$. Combining both estimates, we find:
\[
\tau_c(g, h) = \sup_{\gamma \in \Omega_{C^+}(g, h)} \ell(\gamma) \le \sup_{\gamma \in \Omega_{C^+}(g, h)} L_{\tau_c}(\gamma) \le \sup_{\tilde{\gamma} \in \Omega_{\tau_c}(g, h)} L_{\tau_c}(\tilde{\gamma}) \le \tau_c(g, h),
\]
whence $\tau_c(g, h) = \sup\limits_{\tilde{\gamma} \in \Omega_{\tau_c}(g, h)} L_{\tau_c}(\tilde{\gamma})$.
\end{proof}

\section{Conclusion}

In the present paper, we have obtained a complete characterization of left-invariant intrinsic chronometers on connected finite-dimensional Lie groups, which serves as a sub-Lorentzian analog to V.\,N.~Berestovskii's theorems on left-invariant intrinsic Carnot--Carathéodory--Finsler metrics. We established a one-to-one correspondence between an axiomatized, locally upper semicontinuous, left-invariant intrinsic chronometer $\tau$ and a left-invariant sub-Lorentzian structure $(G, C_0^+, \nu_0)$ defined by a non-trivial closed pointed convex cone $C_0^+ \subset L$ and a concave upper semicontinuous anti-norm $\nu_0$ on $C_0^+$.

In particular, we demonstrated that the condition of compatibility between the chronometer and the Lie group topology---stating that the interiors of the chronological diamonds $\operatorname{int} I(p, q)$ form a basis for the manifold topology---implies the Chow--Rashevskii condition (Hörmander's condition): the linear span of the cone of admissible directions $L_0 = C_0^+ - C_0^+$ generates the entire Lie algebra $L$ ($\widetilde{L}_0 = L$). Furthermore, we established that the introduced chronological condition for the cone $C_0^+$ is a necessary and sufficient condition for the sub-Lorentzian distance functional $\tau_c$, generated by the left-invariant sub-Lorentzian structure $(G, C_0^+, \nu_0)$, to satisfy all the axioms of a left-invariant intrinsic chronometer.

\section*{Acknowledgments}
The author expresses his deep gratitude to his scientific advisor Lev V. Lokutsievskiy for proposing the problem, constant attention to the work, and fruitful discussions that contributed to the improvement of this paper. The author bears sole responsibility for any possible errors and inaccuracies in this manuscript.

\end{document}